\documentclass[courier,notitlepage,12pt]{article}
\usepackage{fancyhdr}
\usepackage{graphics,color}
\usepackage{float}
\usepackage{epsfig}
\usepackage{version}
\usepackage{amsmath}
\usepackage{amssymb}
\usepackage{ntheorem}
\usepackage{amscd}
\usepackage{pst-all}
\usepackage{setspace}
\usepackage{index}
\usepackage{graphicx}

\newtheorem{lemma}{Lemma}[section]
\newtheorem{proposition}{Proposition}[section]
\newtheorem{remark}{Remark}[section]
\newtheorem{theorem}{Theorem}[section]
\newcommand{\preuve}[2]{{\raggedleft \bf Proof#1:} #2  \begin{flushright} \vspace{-8mm} $\blacksquare$ \end{flushright}}
\numberwithin{equation}{section}

\newcommand{\CDD}[3]{\frac{\partial^2 }{\partial #2 \partial #3}#1}

\newcommand{\D}[2]{\frac{\partial }{\partial #2}#1}
\newcommand{\filtr}{\textbf{F}}
\newcommand{\mfiltr}{\textbf{G}}
\newcommand{\st}{\textbf{S}}
\newcommand{\xst}{\textbf{S}}
\newcommand{\expect}{\mathbb E }
\newcommand{\prob}{\mathbb{P}}
\newcommand{\re}{\mathbb{R}}

\newcommand{\open}{\textbf{O}}

\newcommand{\cont}{\textbf{C}}
\newtheorem{assumption}{Assumption}[section]

 \numberwithin{equation}{section}
\allowdisplaybreaks

\usepackage{setspace}

\begin{document}
\bibliographystyle{unsrt}

\title{Impulse Control of a Diffusion with a Change Point}

{\author{ \small
{\bf Lokman A. Abbas-Turki}$^{*}$, {\bf Ioannis Karatzas}$^{\dag}$ and
{\bf Qinghua Li}$^{\ddag}$ \\
\small  $^{*}$TU Berlin, Stochastik und Finanzmathematik,\\ \small Building MA, Str. des 17. Juni 136, 10623 Berlin, Germany\\
\small  $^{\dag}$Department of Mathematics, Columbia University,\\ \small New York, NY 10027, USA\\
\small  $^{\ddag}$Statistics Department, Columbia University,\\ \small 1255 Amsterdam Ave., New York, NY 10027, USA}}

\maketitle

\begin{abstract}
This paper solves a Bayes sequential impulse control problem for a
diffusion, whose drift has an unobservable parameter with a change
point. The partially-observed problem is reformulated into one with
full observations, via  a change of probability measure which
removes the drift. The optimal impulse controls can be expressed in
terms of the solutions and the current values of a Markov process
adapted to the observation filtration. We shall illustrate the
application of our results using the Longstaff-Schwartz algorithm
for multiple optimal stopping times in a geometric Brownian motion
stock price model with drift uncertainty.\bigskip


\end{abstract}

\section{Introduction}\label{sec intr}
\noindent Suppose the price evolution of a stock follows a geometric
Brownian motion, whose drift will change  at an unknown future time
to an unknown level. An investor, who purchases a certain small
number of shares at an initial time, can only observe the evolution
of prices. Based on these observed prices and some reasonable
\textit{a priori} knowledge about the change in the drift, what is
the best time to sell the shares, in order to maximize the
investor's expected profit? Under the same circumstances, what about
the optimal discretely-balanced buying/selling trading strategies
with a larger number of shares?

At a more abstract level, this is a question about optimal stopping
 and impulse control of a diffusion, whose drift term has an unobservable
parameter with a change point. There have been three common
approaches to such problems. The conservative approach is the
mini-max philosophy that optimizes  the worst-case scenario,
formulated as a zero-sum game between a controller and a stopper by
 \cite{KarZam08} Karatzas and Zamfirescu (2008). The approach employed by many  practitioners is to divide model calibration and decision
making into two separate steps. Another approach is to convert the
decision-making problem with partial observations, into one with
full observations, by augmenting the state process with the
posterior probability distribution of the parameter. One
illustration of this method is the work \cite{DZhZh2010} Dai, Zhang
and Zhu (2010). In a geometric Brownian motion price model with
drift uncertainty, the authors assume no impact of the trading
activities, and find two optimal sequences of times to place,
respectively, buy and sell orders.

The topics mentioned in the previous paragraph are all very well
developed fields of research with an extensive literature from the past decades: among them,   \cite{Shiryaev69} Shiryaev (1969) and \cite{Kar03Bayes} Karatzas (2003) for sequential detection;  \cite{Shiryaev78} Shiryaev (1978) or Appendix D in \cite{KarShr98} Karatzas and Shreve (1998) for optimal stopping problems;    \cite{Bens74}, \cite{BensLions74}, \cite{BensLions75a} and \cite{BensLions75b} by Bensoussan and Lions and \cite{OksSul07} {\O}ksendal and Sulem (2007) for impulse controls; as well as  \cite{LiptserShiryaev2001} Liptser and Shiryaev (2001) and \cite{Bens92} Bensoussan (1992) for solving partially observed control problems using filtering techniques.

This paper is an attempt at solving   impulse control problems in
the Bayes sequential framework in one step, without tracking the
posterior probability processes. The conversion from partial to full
observations is facilitated by a change of probability measure,
which hides the drift part of the diffusion. The measure change
method was originally developed for solving change-point detection
problems. For our problem, the state process augmented by the
likelihood ratios is Markovian under the reference probability
measure, and one can derive the dynamic programming principle
satisfied by the value functions. The current values of the
augmented state process provide all the information necessary for
decision making.

There are at least three widely used methods for describing the
value functions of stochastic control problems - PDE, dynamic
programming, and backward SDE. They are different formulations of
the   notion of the ``stochastic maximum principle" (c.f.
\cite{Kushner72} Kushner (1972) and \cite{Davis73} Davis (1973)).
The common mechanism of all  three methods is that the sum of the
value function and the cumulative reward, when evaluated along the
state process corresponding to any admissible control strategy,
yields a supermartingale -- which becomes a martingale if and only
if the control strategy is optimal. After reduction to the Markovian
case, we represent the optimal impulse control via the dynamic
programming principle. Unlike in the pioneering papers \cite{Bens74},
\cite{BensLions74}, \cite{BensLions75a}, \cite{BensLions75b} and
\cite{Bens92} of Bensoussan and Lions, the variational inequalities
associated with the value functions will not be presented here,
because they are numerically inefficient to implement due to the
dimensionality of the augmented state process. From a numerical
point of view, we adapt the Longstaff-Schwartz algorithm for
multiple optimal stopping times. Based on the formulation of the
dynamic programming principle in terms of stopping times, this
method is computationally efficient although it uses Monte Carlo
simulation. We also show the good accuracy of the obtained results
for a simple problem involving geometric Brownian motion and two
optimal stopping times. The reduction to full observation via the
change of measure and the proposed algorithm that involves Monte
Carlo are well suited to a high dimensional state process. They are
a contribution to the methodology of solving partial observation
control problems.

In Section \ref{sec model 1.1}, we specify the model and the quantity that we want
to optimize. Section \ref{sec general theory 2.2} deals with the main theoretical
results related to writing and solving the impulse control problem under a different
probability measure. Section \ref{subsec GBM} illustrates how to use the theoretical
results through the implementation of Longstaff-Schwartz algorithm for multiple optimal
stopping times. Finally, we show in Section \ref{discu} that the method presented
in this paper provides a better framework for multidimensional problems than the usual
one that involves the posterior probabilities.

%%%%%%%%%%%%%%%%%%
\section{Problem formulation}
\label{sec model 1.1}
%%%%%%%%%%%%%%%%%%

This section sets up the diffusion model with a change point in its
drift with Section \ref{subsec 2.1} and lays out the impulse
control problem of the diffusion with Section \ref{subsec 2.2}.

\subsection{The diffusion with a change point}\label{subsec 2.1}

 We consider a probability space $(\Omega,\mathbb{F},\prob^0)$,  which supports
 two independent, $\mathbb{F}$-measurable random variables $\rho$ and $U$, as well as a
 one-dimensional standard  Brownian motion $W^0(\cdot)$ with
respect to its natural filtration $\filtr^{W^0}$. The vector of
random variables $(\rho,  U)$  is independent of the Brownian motion
$W^0(\cdot)$. Let
\begin{equation}\label{larger filtr}
\textbf{G}=\{\textbf{G} (t)\}_{0\leq t\leq T} =\sigma\left\{\rho,
U,W^0(s);0\leq s\leq t\right\}_{0\leq t\leq T}
\end{equation}
denote the filtration generated by $\rho$, $U$ and $W^0(\cdot)$.

The unobservable process
\begin{equation}
\begin{split}
\theta:\; &[0,T]\times\Omega \rightarrow \Theta,\\
 &(t,\omega)\mapsto
\theta(t,\omega)=:\theta(t)
\end{split}
\end{equation}
takes values in the parameter space $\Theta=\{\mu_0,\mu_1,\cdots,\mu_m\}\subset \re$. The process $\theta(\cdot)$ starts with initial value $\theta(0)=\mu_0$, and keeps this value until an unobservable time $\rho$ of {\it regime change}. At the time $\rho$, the parameter $\theta(\cdot)$ changes to a new level $U$, a random variable taking  values in the set $\{ \mu_1, \cdots,\mu_m\}$, and remains at that level until the fixed   terminal time $T\in (0, \infty)$. If   regime change does not occur by time $T$, then $\theta(\cdot)$ takes the value $\mu_0$ throughout the interval $[0,T]$, namely,
\begin{equation} \label{theta 1.1}
\theta(t)=\left\{ \begin{aligned}
         \mu_0&,\;0\leq t<\rho\wedge T\,; \\
         U&,\;\rho\wedge T\leq t\leq T\,.
                          \end{aligned} \right.
                          \end{equation}
The change point $\rho$ and the level $U$ have prior distributions
\begin{equation}
\label{rho prior 1.1}
\prob(\rho>t)=e^{-\lambda t},\;t\geq 0,
\end{equation}
 and
\begin{equation} \label{theta prior 1.1}
\prob(U=\mu_j)=p_j,\;j=1,2,\cdots,m.
\end{equation}
For any possible values $u\in\Theta$, the given measurable functions
$b(\cdot,\cdot\,;u):[0,T]\times \re \rightarrow \re$ and
$\sigma(\cdot,\cdot):[0,T]\times \re \rightarrow \re$    satisfy the
following Lipschitz and boundedness condition.

\begin{assumption}
\label{assump sde 1.1}
There exists a constant $C>0$, such that\\
(i) for all $(t,x^1)$, $(t,x^2)\in [0,T]\times \re$, and for all $u \in \Theta$, we have
\begin{eqnarray}\label{cond Lip}
\begin{array}{ccc} \displaystyle
\hspace{-4cm}|b(t,x^1;u)-b(t,x^2;u)|+|\sigma(t,x^1)-\sigma(t,x^2)| \vspace{1mm}\\
\hspace{3.1cm}+\displaystyle\left|\frac{b(t,x^1;u)}{\sigma(t,x^1)}-\frac{b(t,x^2;u)}{\sigma(t,x^2)}\right|
\leq C|x^1-x^2|;
\end{array}
\end{eqnarray}
whereas \\
(ii) for all $(t,x)\in [0,T]\times \re$ and  all $u \in \Theta$, we have
\begin{equation}\label{cond bdd}
\sigma(t,x)>0 \qquad \text{ and } \qquad \left|\frac{b(t,x;u)}{\sigma(t,x)}\right|\leq C\,.
\end{equation}
\end{assumption}

The Assumption \ref{assump sde 1.1} (i)  implies a linear growth
condition on the functions $b(t,\cdot\,;u)$ and
$\sigma(t,\cdot\,)\,$: there exists another constant $C'>0$, such
that
\begin{equation}
\label{cond linear}
\left|b(t,x;u)\right|+\left|\sigma(t,x)\right|\leq C'|1+x|\;
\end{equation}
holds for all $u \in \Theta$ and all $(t,x)\in [0,T]\times \re$.

\smallskip
Let $N$ be a positive integer, $0=\tau_0< \tau_1<\tau_2<\cdots
<\tau_N \leq T$ be stopping times with respect to the filtration
$\filtr^{W^0}\,$, and $\zeta_i$ be an $\re$-valued $\filtr^{W^0}(\tau_i-)$
- measurable random variable for $i=1,2,\cdots,N$. The $N$-tuple
$(\tau,\zeta)=\{(\tau_i,\zeta_i)\}_{i=1}^N$ is called an {\it
impulse control}. The set of {\it admissible controls}, denoted as
$\textbf{I}$, is the collection of all
 such impulse controls $(\tau,\zeta)$.
The jump size $\gamma:\re\times\re\rightarrow\re$ is a given bounded,
measurable function. Given an arbitrary impulse control $(\tau,\zeta) \in \textbf{I}$,
the controlled state process $X(\cdot)$ is the unique strong solution to the equation
\begin{equation}\label{sde 1.2}
X(t)=x_0+\int_0^t\sigma(s,X(s))dW^0(s) + \sum\limits_{\tau_i\leq
t}\gamma(X(\tau_i-),\zeta_i),\;\quad 0\leq t\leq T.
\end{equation}
By Assumption \ref{assump sde 1.1} (ii)  and equation (\ref{sde
1.2}), the Brownian filtration $\filtr^{W^0}$ coincides with
$\filtr=\{\filtr (t)\}_{0\leq t\leq T}$, which denotes the
filtration generated by the process $X(\cdot)$. The collection of all
$\filtr$-stopping times with values in $[0,T]$ is denoted as $\xst$, and the collection of all
$\filtr$-stopping times with values in $[t,T]$ is denoted as $\xst_t$.

\medskip
\noindent $\bullet~$ In order to define another probability measure
$\prob$ on the space $\Omega$ and on the sigma algebra
$\textbf{G}(T)\,$, the one with respect to which we shall formulate
our impulse control problem, we introduce  the $\textbf{G}$-adapted
process
\begin{equation}\label{RN deriv 2.1}
Z(t)=\exp \left\{\int_0^t
\frac{b(s,X(s);\theta(s))}{\sigma(s,X(s))}dW^0(s)-\frac{1}{2}\int_0^t
\frac{b^2(s,X(s);\theta(s))}{\sigma^2(s,X(s))}ds\right\}\text{,
}0\leq t\leq T\,,
\end{equation}
which will play the role of Radon-Nikodym derivative of the new
measure $\prob$ with respect to the ``reference probability measure"
$\prob^0$. Whereas, for every number $u \in \Theta$,  the
$\filtr$-adapted likelihood ratio process is defined as
\begin{equation}\label{lr 2.1}
L(t;u)=\exp \left\{\int_0^t
\frac{b(s,X(s);u)}{\sigma(s,X(s))}dW^0(s)-\frac{1}{2}\int_0^t
\frac{b^2(s,X(s);u)}{\sigma^2(s,X(s))}ds\right\},\;0\leq t\leq T\,.
\end{equation}
From the expression (\ref{theta 1.1}) for $\theta(\cdot)$, the
Radon-Nikodym derivative $Z(\cdot)$ can be written, in terms of the
likelihood ratio process $L(\cdot;u)$ and of the random vector
$(\rho,U)$, as
\begin{equation}
\label{RN deriv 1.2}
Z(t)\,=\,L(\rho;\mu_0)\left(\sum_{j=1}^m
\mathbf{ 1}_{\{U=\mu_j\}}\frac{L(t;\mu_j)}{L(\rho;\mu_j)}\right)\mathbf{ 1}_{\{\rho<t\}}
+L(t;\mu_0)\mathbf{ 1}_{\{\rho\geq t\}} \,,\;~~~~0\leq t\leq T.
\end{equation}
The Radon-Nikodym process $Z(\cdot)$ in (2.10) is a
$(\prob^0,\textbf{G})$-martingale, because of Assumption \ref{assump sde 1.1} (ii)  on the boundedness of the ratio
$b(\cdot\,,\cdot\,;u)/\sigma(\cdot,\cdot)$ and of the Novikov condition; the same is true for  the likelihood ratio process $L(\cdot\,;u)$ in (2.11), for any $u\in \Theta$. There exists then a probability measure $\prob$ equivalent to $\prob^0$, satisfying
\begin{equation}
\label{change measure 1.1}
\left.\frac{d\prob}{d\prob^0}\right|_{\mfiltr(t)}=Z(t),\;~~~~0\leq
t\leq T.
\end{equation}
Under this new probability measure $\prob$, the random variables
$\rho$ and $U$ are still independent and retain the prior
distributions of (\ref{rho prior 1.1}) and (\ref{theta prior 1.1}).
By a generalization of the Girsanov theorem to local martingales in
\cite{SchuppenWong74} Van Schuppen and Wong (1974), the process
\begin{equation}
\left\{\int_0^t\sigma(s,X(s))dW^0(s)-\int_0^t b(s,X(s);\theta(s))ds\right\}_{0\leq t\leq T}
\end{equation}
is a local $(\prob,\mfiltr)$-martingale, having the instantaneous
quadratic variation $\sigma^2(\cdot\,,X(\cdot))$; and the process
$W(\cdot)$ defined as
\begin{equation}
W(t):=W^0(t)-\int_0^t
\sigma^{-1}(s,X(s))b(s,X(s);\theta(s))ds\,,\;~~~~0\leq t\leq T
\end{equation}
is a standard $\prob$-Brownian
motion. The process $X(\cdot)$ defined by (\ref{sde 1.2})
is also the unique strong solution to the equation
\begin{equation}\label{SDE impulse 2.1}
X(t)=x_0+\int_0^t \hspace{-2mm} b(s,X(s);\theta(s))ds+\int_0^t\hspace{-2mm}
\sigma(s,X(s))dW(s) + \sum \limits_{\tau_i\leq
t}\hspace{-1mm}\gamma(X(\tau_i-),\zeta_i),\;0\leq t\leq T.
\end{equation}

\subsection{Impulse control of the diffusion}\label{subsec 2.2}

The impulse control problem we study in this paper, consists in
choosing an optimal impulse control
$(\tau^*,\zeta^*)=\{(\tau_i^*,\zeta_i^*)\}_{i=1}^N$  to achieve the
maximal expected reward
\begin{equation}
\label{opt value 2.1.0}
\begin{split}
V:=\sup\limits_{(\tau,\zeta)\in\textbf{I}}\expect\left[\int_0^T
h(X(t))dt+\xi(X(T))+\sum\limits_{i=1}^N
c(X(\tau_i-),\zeta_i)\right],
\end{split}
\end{equation}
over all admissible impulse controls
$(\tau,\zeta)=\{(\tau_i,\zeta_i)\}_{i=1}^N$ in $\textbf{I}$. The reward
functions $\xi$ and $h:\re\rightarrow \re$ are measurable and
satisfy  conditions  (i)  and  (ii)  in Assumption \ref{assump rwd
1.1} below. Furthermore, we impose growth condition on the
deterministic measurable functions $\gamma$ and
$c:\re\times\re\rightarrow\re$ in the state variable, as follows.
\begin{assumption}
\label{assump rwd 1.1} (i) The function $\xi(\cdot)$ is twice
continuously differentiable, with first and second order derivatives
denoted as $\xi'(\cdot)$ and
$\xi''(\cdot)$.\\
(ii) The functions $h(\cdot)$, $\xi(\cdot)$, $\xi'(\cdot)$ and
$\xi''(\cdot)$ are locally Lipschitz and have polynomial growth.\\
(iii) The function $\gamma(x,z)$ is bounded for all $x\in\re$ and
$z\in\re$,
 and the function $c(x,z)$ has
polynomial growth rate in $x\in\re$ uniformly for all $z\in\re$.
Both functions $\gamma(x,z)$ and $c(x,z)$ are continuous in $x$, for
any arbitrarily fixed $z\in\re$.
\end{assumption}

%%%%%%%%%%%%%%%%%%%%%%
\section{Solving the impulse control problem}
\label{sec general theory 2.2}
%%%%%%%%%%%%%%%%%%%%%%%

This section provides a theoretical solution to the impulse control
problem (\ref{opt value 2.1.0}): Section \ref{subsec measure change}
reduces the partially observable problem into one of full
observation, by changing to the reference probability measure under
which the state process is a martingale and augmenting the state
process with the likelihood ratios and their integrals. Section
\ref{subsec solution} solves the fully observable impulse control
problem under the reference probability measure, by representing the
optimal control in terms of the value functions and the augmented
state process.

\subsection{The measure change method}\label{subsec measure change}

By the ``measure change method", we mean considering the impulse
control problem (\ref{opt value 2.1.0}) under the reference
probability measure $\prob^0$ which removes the unobservable drift
of the state process $X(\cdot)$. Using the Bayes rule and the
properties of conditional expectations, the maximal expected reward
$V$ from (\ref{opt value 2.1.0}) can be written as
\begin{equation}
\label{opt value 2.1}
\begin{split}
V=&\sup\limits_{(\tau,\zeta)\in\textbf{I}}\expect^0\left[Z(T)\left(\int_0^T
h(X(t))dt+\xi(X(T))+\sum\limits_{i=1}^N c(X(\tau_i-),\zeta_i)\right)\right]\\
=&\sup\limits_{(\tau,\zeta)\in\textbf{I}}\expect^0\left[
\expect^0\left[Z(T)\left|\filtr(T)\right.\right] \left(\int_0^T
h(X(t))dt+\xi(X(T))+\sum\limits_{i=1}^N
c(X(\tau_i-),\zeta_i)\right)\right].
\end{split}
\end{equation}
From the Bayes point of view, the quantity
$\expect^0\left[Z(t)\left|\filtr (t)\right.\right]$ in (\ref{opt
value 2.1}) is the posterior expectation of the Radon-Nikodym
derivative $Z(\cdot)$ under the reference probability measure
$\prob^0$, given the observations of $X(\cdot)$ up to date. Because
of the independence of $(\rho,\,U)$ and $X(\cdot)$ under $\prob^0$,
from the prior $\prob^0$-distributions (\ref{rho prior 1.1}) and
(\ref{theta prior 1.1}), and by (\ref{RN deriv 1.2}), this posterior
expectation has the form
\begin{equation}\label{post expect 1.1}
\begin{split}
\expect^0\left[Z(t)\left|\filtr (t)\right.\right] =\sum_{j=1}^m
\left(p_j L(t;\mu_j)\int_0^t\frac{L(s;\mu_0)}{L(s;\mu_j)} \lambda
e^{-\lambda s}ds\right) +e^{-\lambda t}L(t;\mu_0) ,\;0\leq t\leq
T.
\end{split}
\end{equation}
For every $u \in \Theta$, the likelihood ratio process
$L(\cdot\,;u)$ defined in (\ref{lr 2.1}) is a
$(\prob^0,\filtr)$-martingale
satisfying the stochastic integral equation
\begin{equation}\label{dlr 2.1}
\begin{split}
L(t;u)=\int_0^t L(s;u)\frac{b(s,X(s);u)}
{\sigma(s,X(s))}dW^0(s)\,\text{, }~~~~0\leq t\leq T,
\end{split}
\end{equation}
with respect to the standard $(\prob^0,\filtr)$-Brownian motion
$W^0(\cdot)$. From equations (\ref{post expect 1.1}) and (\ref{dlr
2.1}) we obtain, for $0\leq t\leq T$, that
\begin{eqnarray}\label{post expect 1.2}
\hspace{-2mm}\begin{array}{ccc}
d\left(\expect^0\left[Z(t)\left|\filtr (t)\right.\right]\right)
\hspace{-2mm} &=& \hspace{-3mm}\displaystyle \hspace{-1mm}\sum_{j=1}^m p_j\left(\int_0^t\frac{L(s;\mu_0)}{L(s;\mu_j)}
\lambda e^{-\lambda s}ds\right)dL(t;\mu_j) +e^{-\lambda
t}dL(t;\mu_0) \vspace{2mm} \\
\hspace{-2mm} &= & \hspace{-3mm}\left(\hspace{-2mm}\begin{array}{c}
\displaystyle \sum_{j=1}^m p_j\left(\int_0^t \hspace{-1mm}\frac{L(s;\mu_0)}{L(s;\mu_j)}
\lambda e^{-\lambda s}ds\right)L(t;\mu_j)\frac{b(t,X(t);\mu_j)}
{\sigma(t,X(t))} \\ \displaystyle + e^{-\lambda t}L(t;\mu_0)\frac{b(t,X(t);\mu_0)}
{\sigma(t,X(t))}\end{array}\hspace{-2mm}\right)\hspace{-1mm}dW^0(t),
\end{array}\hspace{-5mm}
\end{eqnarray}
so the posterior expectation $\left\{\expect^0\left[Z(t)\left|\filtr
(t)\right.\right]\right\}_{0\leq t\leq T}$ is a local
$(\prob^0,\filtr)$-martingale; in fact a
$(\prob^0,\filtr)$-martingale, as is easily checked from the
definition (2.13).\\

Applying It\^{o}'s formula, we shall see that the random variable
inside the $\prob^0$-expectation on the second line of (\ref{opt
value 2.1}) is a $(\prob^0,\filtr)$-semimartingale evaluated at the
time $T$. Lemma \ref{lem mart 1.1} will show that its local
martingale part is in fact a square-integrable
$(\prob^0,\filtr)$-martingale of class $\textbf{D}$ (Definition 4.8,
page 24 of \cite{KarShr88} Karatzas and Shreve (1988)), because of
the uniform integrability result from Lemma \ref{lem UI 1.1}.

\begin{lemma}
\label{lem UI 1.1}
For every $u\in\Theta$, consider the process $R(\cdot\,;u)$  defined as
\begin{equation}\label{lr 1.2}
R(t;u):=\int_0^t\frac{L(s;\mu_0)}{L(s;u)} \lambda
e^{-\lambda s}ds\text{, }~~~~0\leq t\leq T.
\end{equation}
For any nonnegative integers $q_1$, $q_2$ and $q_3$, we have
\begin{equation}\label{UI est 1.1}
\expect^0\left[
\sup\limits_{0\leq t\leq
T}|X(t)|^{q_1}\right]<\infty\text{, }
\expect^0\left[
\sup\limits_{0\leq t\leq
T}L^{q_2}(t;u)\right]<\infty~~\text{ and }~~
\expect^0\left[R^{q_3}(T;u)\right]<\infty\,;
\end{equation}
furthermore, the family
\begin{equation}\label{UI 1.2}
\left\{\sup\limits_{0\leq
t\leq T}|X(t)|^{q_1}L^{q_2}(\tau;u)R^{q_3}(\tau;u)\right\}_{\tau\in\xst}
\end{equation}
is uniformly integrable with respect to the probability measure
$\prob^0$.
\end{lemma}
\preuve{}{By Assumption \ref{assump rwd 1.1} (iii)  and
equation (\ref{sde 1.2}), there exists $
C_1(x_0,\gamma,N,q_1)\in (0, \infty)$, such that for any $0\leq t
\leq T$, we have
\begin{equation}\label{est 1.12}
\begin{split}
&\sup\limits_{0\leq s\leq t}|X(s)|^{2{q_1}}
\leq  \sup\limits_{0\leq s\leq t}\left(x_0+\left|\int_0^s\sigma(w,X(w))dW^0(w)\right|
+ \sum\limits_{i=1}^N \left|\gamma(X(\tau_i-),\zeta_i)\right|\right)^{2{q_1}}\\
\leq & C_1(x_0,\gamma,N,q_1)\left(1+
\sup\limits_{0\leq s\leq t}\left|\int_0^s\sigma(w,X(w))dW^0(w)\right|^{2{q_1}}\right).
\end{split}
\end{equation}
Since $\int_0^\cdot\sigma(t,X(t))dW^0(t)$ is a local $\prob^0$-martingale, from the Burkholder-Davis-Gundy inequality (e.g. page 166 of \cite{KarShr88}
Karatzas and Shreve (1988)),  for $q_1=1,2,\cdots$ we have
\begin{equation}\label{est 1.5}
\expect^0\left[\sup\limits_{0\leq t\leq T}\left|\int_0^t\sigma(s,X(s))dW^0(s)\right|^{2{q_1}}\right]\leq
C_2(q_1)\expect^0\left[\left(\int_0^T\sigma^2(t,X(t))dt\right)^{q_1}\right],
\end{equation}
for some constant $0< C_2(q_1)<\infty$. But there exists a constant $ C_3(\sigma,q_1,T)\in (0, \infty)$, such that
\begin{equation}\label{est 1.7}
\begin{split}
&\expect^0\left[\left(\int_0^T\sigma^2(t,X(t))dt\right)^{q_1}\right]
\leq  T^{q_1-1}\expect^0\left[\int_0^T\sigma^{2q_1}(t,X(t))dt\right]\\
= &T^{q_1-1}\int_0^T\expect^0\left[\sigma^{2q_1}(t,X(t))\right]dt
\leq T^{q_1-1}C_3(\sigma,q_1,T)\left(1+\int_0^T\expect^0\left[|X(t)|^{2q_1}\right]dt\right)\\
\leq &
T^{q_1-1}C_3(\sigma,q_1,T)\left(1+\int_0^T\expect^0\left[\sup\limits_{0\leq
s\leq t}|X(s)|^{2{q_1}}\right]dt\right),
\end{split}
\end{equation}
where the second inequality comes from the inequality (\ref{cond linear}), the linear growth property of
$\sigma(t,\cdot)$. Inequalities (\ref{est
1.12}), (\ref{est
1.5}) and (\ref{est 1.7}) imply
\begin{equation}\label{est 1.8}
\begin{split}
&\expect^0\left[\sup\limits_{0\leq t\leq T}\left|\int_0^t\sigma(s,X(s))dW^0(s)\right|^{2{q_1}}\right]\\
\leq \, &\, C_4(x_0,\sigma,\gamma,N,q_1,T)\left(1+\int_0^T\expect^0\left[\sup\limits_{0\leq
s\leq t}\left|\int_0^s\sigma(v,X(v))dW^0(v)\right|^{2{q_1}}\right]dt\right),
\end{split}
\end{equation}
for some constant $ C_4(x_0,\sigma,\gamma,N,q_1,T)\in (0, \infty)$.
Then, by the Gronwall inequality (e.g. page 287 of \cite{KarShr88}
Karatzas and Shreve (1988)), we know that
\begin{equation}
\expect^0\left[\sup\limits_{0\leq t\leq T}\left|\int_0^t\sigma(s,X(s))dW^0(s)\right|^{2{q_1}}\right]<\infty,
\end{equation}
hence by inequality (\ref{est
1.12}) we have
\begin{equation}\label{est 1.9}
\expect^0\left[\sup\limits_{0\leq t\leq T}|X(t)|^{2{q_1}}\right]<
\infty.
\end{equation}
The equation (\ref{lr 2.1}) and Assumption \ref{assump sde 1.1}
(ii) imply that, for $0\leq t\leq T$,
\begin{equation}
\label{dlr 2.3}
\tilde{L}(t;u)\leq L^{-1}(t;u)=\exp \left\{ \int_0^t
\frac{b^2(s,X(s);u)}{\sigma^2(s,X(s))}ds\right\}\tilde{L}(t;u)\leq
\exp\{C^2 T\}\tilde{L}(t;u),
\end{equation}
where we have defined
\begin{equation}\label{dlr 2.4}
\tilde{L}(t;u):=\exp \left\{-\int_0^t
\frac{b(s,X(s);u)}{\sigma(s,X(s))}dW^0(s)-\frac{1}{2}\int_0^t
\frac{b^2(s,X(s);u)}{\sigma^2(s,X(s))}ds\right\}
\end{equation}
and noted
\begin{equation}\label{dlr 2.2}
\tilde{L}(t;u)=-\int_0^t \tilde{L}(s;u)\frac{b(s,X(s);u)}
{\sigma(s,X(s))}dW^0(s).
\end{equation}
Using Assumption \ref{assump sde 1.1} (ii)  and the same arguments
as those leading to (\ref{est 1.9}), as well as the equations
(\ref{dlr 2.1}), (\ref{dlr 2.3}) and (\ref{dlr 2.2}), we can show
\begin{equation}
\label{est 1.10}
\expect^0\left[\sup\limits_{0\leq t\leq T}L^{q_2}(t;u)\right]<
\infty~~~\text{ and }~~~\expect^0\left[\sup\limits_{0\leq t\leq
T}L^{-q_2}(t;u)\right]<\infty.
\end{equation}
It follows from equations (\ref{lr 1.2}) and (\ref{est 1.10}) that
\begin{equation}
\label{est 1.11}
\expect^0\left[R^{q_3}(T;u)\right]<
\infty.
\end{equation}
Taking an arbitrary
$\filtr$-stopping time $\tau$ with values in $[0,T]$, by the H\"older
inequality and by the estimates
(\ref{est 1.9}), (\ref{est 1.10}) and (\ref{est 1.11}), we get
\begin{equation}
\label{est 1.2}
\begin{split}
&\expect^0\left[\sup\limits_{0\leq
t\leq T}|X(t)|^{q_1} L^{q_2}(\tau;u)R^{q_3}(\tau;u)  \right]\\
\leq &\left(\expect^0\left[\sup\limits_{0\leq t\leq
T}|X(t)|^{2{q_1}} \right] \right)^{1/2}\left(\expect^0\left[L^{4q_2}(\tau;u)\right]
\right)^{1/4}\left(\expect^0\left[R^{4{q_3}}(\tau;u)\right]
\right)^{1/4}\\
\leq &\left(\expect^0\left[\sup\limits_{0\leq t\leq
T}|X(t)|^{2{q_1}} \right] \right)^{1/2}\left(\expect^0\left[\sup\limits_{0\leq t\leq
T}L^{4q_2}(t;u)\right]
\right)^{1/4}\left(\expect^0\left[R^{4{q_3}}(T;u)\right]
\right)^{1/4}<\infty.
\end{split}
\end{equation}
To derive the uniform integrability of the family (\ref{UI 1.2})
from (\ref{est 1.2}), we use the Cauchy-Schwartz and Chebyshev
inequalities to get the estimate
\begin{equation}
\begin{split}
&\sup\limits_{\tau\in\xst}\, \expect^0\left[
\sup\limits_{0\leq
t\leq T}|X(t)|^{q_1} L^{q_2}(\tau;u)R^{q_3}(\tau;u)
\mathbf{ 1}_{ \big\{ \sup\limits_{0\leq
t\leq T}|X(t)|^{q_1} L^{q_2}(\tau;u)R^{q_3}(\tau;u)>A \big\}}\right]\\
\leq &  \sup\limits_{\tau\in\xst}\left(\expect^0\left[
\sup\limits_{0\leq t\leq T}|X(t)|^{2q_1}
L^{2q_2}(\tau;u)R^{2q_3}(\tau;u)\right]\right)^{1/2}\cdot\\
&\cdot\left(\prob^0\left( \sup\limits_{0\leq t\leq T}|X(t)|^{q_1}
L^{q_2}(\tau;u)R^{q_3}(\tau;u)>A\right)\right)^{1/2}\\
\leq  &\frac{1}{A}\sup\limits_{\tau\in\xst}\expect^0\left[
\sup\limits_{0\leq t\leq T}|X(t)|^{2q_1}
L^{2q_2}(\tau;u)R^{2q_3}(\tau;u)\right],
\end{split}
\end{equation}
which tends to zero as $A\rightarrow\infty$, on the strength of (\ref{est 1.2}).}

\begin{lemma}\label{lem mart 1.1} For $0\leq t\leq T$, $x\in\re$,
$l=(l_0,l_1,\cdots,l_m)\in \re^{m+1}$, and $r=(r_1,\cdots,r_m)\in
\re^m$, consider the function $\alpha$   defined as
\begin{equation}\label{cross var 1.1}
\begin{split}
\alpha(t,x,l,r) :=&\left(\sum\limits_{j=1}^m p_jl_jr_j+e^{-\lambda
t}l_0\right)\left(h(x)+\frac{1}{2}\xi''(x)\sigma^2(t,x)\right)
\\
&+\left(\sum\limits_{j=1}^m p_jl_jr_jb(t,x;\mu_j)+e^{-\lambda
t}l_0b(t,x;\mu_0)\right)\xi'(x),
\end{split}
\end{equation}
and the function $\beta$  defined as
\begin{equation}\label{cross var 2.1}
\begin{split}
\beta(t,x,l,r,z) := \hspace{-1mm} \left( \hspace{-1mm} \sum\limits_{j=1}^m
p_jl_jr_j+e^{-\lambda t}l_0 \hspace{-1mm} \right)\hspace{-1mm}
\big(\xi(x+\gamma(x,z))-\xi(x)+\xi'(x)\gamma(x,z)+c(x,z)\big)\,.
\end{split}
\end{equation}
Then, for $0\leq t\leq T$, we have
\begin{equation}\label{Dyn 1.0}
\begin{split}
&\expect^0\left[Z(t)\left|\filtr (t)\right.\right]\left(\int_0^t
h(X(s))ds+\xi(X(t))+\sum\limits_{\tau_i\leq
t}c(X(\tau_i-),\zeta_i))\right)\\
=\,&M^0(t)+\int_0^t\alpha \big(s,X(s),L(s),
R(s)\big)ds+\sum\limits_{\tau_i\leq
t}\beta \big(\tau_i,X(\tau_i-),L(\tau_i-),R(\tau_i-),\zeta_i\big),
\end{split}
\end{equation}
where $M^0(\cdot)$ is some square integrable
$(\prob^0,\filtr)$-martingale with $M^0(0)=\xi(x_0)$,
\begin{equation}\label{lr 1.3}
R(t):=\left(R(t;\mu_1)\cdots,R(t;\mu_m)\right),
\end{equation}
and
\begin{equation}\label{lr 1.4}
L(t):=\left(L(t;\mu_0),L(t;\mu_1)\cdots,L(t;\mu_m)\right).
\end{equation}
\end{lemma}
\preuve{}{ Applying It\^{o}'s formula for semimartingales with
jumps, we get
\begin{equation}\label{Dyn 1.1}
\begin{split}
&\expect^0\left[Z(t)\left|\filtr (t)\right.\right] \left(\int_0^t
h(X(s))ds+\xi(X(t))+\sum\limits_{\tau_i\leq
t}c(X(\tau_i-),\zeta_i)\right)
\\
=& \hspace{1mm} \xi(x_0)+\int_0^t\left(\int_0^{s-} \hspace{-2mm}
h(X(u))du+\xi(X(s-))+\hspace{-2mm}\sum\limits_{\tau_i\leq
s-}c(X(\tau_i-),\zeta_i)\right)d
\expect^0\left[Z(s)\left|\filtr (s)\right.\right]\\
&+\int_0^t\expect^0\left[Z(s-)\left|\filtr (s-)\right.\right]\xi'(X(s-))\sigma(t,X(t))dW^0(t)\\
&+\int_{0+}^t \alpha\big(s-,X(s-),L(s-),
R(s-)\big)ds+ \sum\limits_{\tau_i\leq
t}\beta \big(\tau_i,X(\tau_i-),L(\tau_i-),R(\tau_i-),\zeta_i\big).
\end{split}
\end{equation}
By change of variables and the continuity of Riemann integrals,
\begin{equation}\label{Dyn 1.2}
\int_{0+}^t \alpha\left(s-,X(s-),L(s-), R(s-)\right)ds
\,=\,\int_{0}^t \alpha\left(s,X(s),L(s), R(s)\right)ds.
\end{equation}
Define
\begin{equation}\label{mart 2.1}
\hspace{-1mm}\begin{split}
M^0(t) := \xi(x_0)+ \hspace{-1mm} &\int_0^t\hspace{-1mm}\left(\hspace{-1mm}
\int_0^{s-} \hspace{-3mm} h(X(u))du+\xi(X(s-))+ \hspace{-2mm}
\sum\limits_{\tau_i\leq s-}\hspace{-2mm}c(X(\tau_i-),\zeta_i) \hspace{-1mm} \right)
\hspace{-1mm}d \expect^0\left[Z(s)\left|\filtr (s)\right.\right] \\
+&\int_0^t\expect^0\left[Z(s-)\left|\filtr
(s-)\right.\right]\xi'(X(s-)) \sigma(t,X(t))dW^0(t).
\end{split}
\end{equation}
Equations (\ref{Dyn 1.1})-(\ref{mart 2.1}) imply that  (\ref{Dyn 1.0}) holds.
Substituting  (\ref{post expect 1.2}) into   (\ref{mart 2.1}), we get
\begin{eqnarray*}
\begin{array}{c}
\hspace{-3.8cm} M^0(t)\,=  \displaystyle \xi(x_0)+ \int_0^t dW^0(s)
\Bigg[\expect^0\left[Z(s-)\left|\filtr (s-)\right.\right]\xi'(X(s-))
\vspace{2mm} \\ \hspace{-2mm} + \left(\hspace{-1mm}\begin{array}{c}
\displaystyle \sum_{j=1}^m p_j\left(\int_0^t
\hspace{-1mm}\frac{L(s;\mu_0)}{L(s;\mu_j)} \lambda e^{-\lambda
s}ds\right)L(t;\mu_j)\frac{b(t,X(t);\mu_j)} {\sigma(t,X(t))} \\
\displaystyle + e^{-\lambda t}L(t;\mu_0)\frac{b(t,X(t);\mu_0)}
{\sigma(t,X(t))}\end{array}\hspace{-1mm}\right) \hspace{-2mm} \left(
\hspace{-1mm}\begin{array}{c}\int_0^{s-} h(X(u))du\\
\displaystyle+\xi(X(s-))\\ \displaystyle + \sum\limits_{\tau_i\leq
s-} \hspace{-2mm} c(X(\tau_i-),\zeta_i)\end{array} \hspace{-1mm}
\right)\Bigg].
\end{array}
\end{eqnarray*}
By Lemma \ref{lem
UI 1.1}, $M^0(\cdot)$ is an integral of $\prob^0$-square integrable processes with respect to the $(\prob^0,\filtr)$-Brownian motion $W^0(\cdot)$, hence $M^0(\cdot)$ is also a local $(\prob^0,\filtr)$-martingale.

We need to show that $M^0(\cdot)$ is a $(\prob^0,\filtr)$-martingale,   not just a local martingale. It suffices to show that the family $\{M^0(\tau)\}_{\tau\in\xst}$ is uniformly integrable under the probability measure $\prob^0$. By equations (\ref{post expect 1.1}) and (\ref{Dyn 1.0}), $M^0(\cdot)$ can be expressed alternatively as
$$
M^0(t)\,=\, \left(\sum_{j=1}^m \left(p_j
L(t;\mu_j)\int_0^t\frac{L(s;\mu_0)}{L(s;\mu_j)} \lambda e^{-\lambda
s}ds\right) +e^{-\lambda
t}L(t;\mu_0)\right)\left( \begin{array}{c} \int_0^t h(X(s))ds\\+\xi(X(t)) \end{array}\right)
$$
\vspace{-3mm}
\begin{equation}\label{mart 1.2}
~~~~~~~~~~~~~-\int_0^t \alpha\big(s,X(s),L(s), R(s)\big)ds-\sum\limits_{\tau_i\leq
t}\beta \big(\tau_i,X(\tau_i-),L(\tau_i-),R(\tau_i-),\zeta_i\big)\,.
\end{equation}
From the expressions (\ref{mart 1.2}), (\ref{cross var 1.1}),
(\ref{cross var 2.1}), (\ref{cond linear}) and Assumption
\ref{assump rwd 1.1} (ii)(iii), we know that there exist
a constant $C>0$ and a positive integer $q$, such that
\begin{equation}
\begin{split}
&\big|M^0(t)\big|
\\
\leq &
  C  \left( \begin{array}{c}\sum_{j=1}^m L(t;\mu_j)R(t;\mu_j)+L(t;\mu_0) \vspace{1mm}\\
  +\int_0^t \sum_{j=1}^m L(s;\mu_j)R(s;\mu_j)+L(s;\mu_0)\,ds \end{array} \right)\sup\limits_{0\leq s\leq
T}|X(s)|^q,
\end{split}
\end{equation}
for all $(t,\omega)\in [0,T]\times\Omega$. Then, from Lemma \ref{lem
UI 1.1}, we know that, under the probability measure $\prob^0$, the
local martingale $M^0(\cdot)$ is both square-integrable and of class
$\textbf{D}$ on $[0,T]$. The latter implies that $M^0(\cdot)$ is a
$(\prob^0,\filtr)$-martingale. }

\noindent Because the process $M^0(\cdot)$ in (\ref{Dyn 1.0}) and
(\ref{mart 2.1}) is a $(\prob^0,\filtr)$-martingale, it should
vanish from inside the $\prob^0$-expectation of (\ref{opt value
2.1}), leaving only the initial value and the finite variation part
of the semimartingale. This property enables Lemma \ref{lem opt
value 2.1} to rewrite the $\prob^0$-expectations in (\ref{opt value
2.1}) in a more convenient manner.
\begin{lemma}\label{lem opt value 2.1} For any impulse control $(\tau,\zeta)\in\textbf{I}$,
\begin{equation}
\label{opt value 2.3} \expect^0\left[
\expect^0\left[Z(T)\left|\filtr(T)\right.\right] \left(\int_0^T
h(X(t))dt+\xi(X(T))+\sum\limits_{i=1}^N
c(X(\tau_i-),\zeta_i)\right)\right]
\end{equation}
$$
=\,\xi(x_0)+\expect^0\left[\int_0^T \alpha\left(s,X(s),L(s),
R(s)\right)ds+\sum\limits_{i=1}^N
\beta(\tau_i,X(\tau_i-),L(\tau_i-),R(\tau_i-),\zeta_i)\right].
$$
\end{lemma}
\preuve{}{This is because the process
\begin{equation}
\begin{split}
M^0(t)=&~\expect^0\left[Z(t)\left|\filtr
(t)\right.\right]\left(\int_0^t h(X(s))ds+\xi(X(t))\right)-\int_0^t
\alpha\left(s,X(s),L(s),
R(s)\right)ds\\
&-\sum\limits_{\tau_i\leq
t}\beta(\tau_i,X(\tau_i-),L(\tau_i-),R(\tau_i-),\zeta_i)\,\text{, }\quad
0\leq t\leq T
\end{split}
\end{equation} is a
$(\prob^0,\filtr)$-martingale by Lemma \ref{lem mart 1.1}, hence
$\expect^0\left[M^0(T)\right]=M^0(0)=\xi(x_0)$. Equivalently,
equation (\ref{opt value 2.3}) holds. }

\noindent Up to this point, the $\prob^0$-expected reward from (\ref{opt value
2.1}) has been rewritten into the $\prob^0$-expectation of the sum
of a reward $\alpha$ cumulated over the time interval $[0,T]$, and
of a reward $\beta$ received only at the times of intervention, as
in equation (\ref{opt value 2.3}). Both $\alpha$ and $\beta$ are
functions of the processes $X(\cdot)$, $L(\cdot)$ and $R(\cdot)$,
which are adapted to the observation filtration $\filtr$. Lemma
\ref{lem Markov 2.1} and Proposition \ref{prop non explosion} will
show that the triple of processes $(X(\cdot),L(\cdot),R(\cdot))$
forms a well-behaved Markov process, because it is the unique strong
solution to a stochastic differential equation with locally
Lipschitz coefficients and this solution does not explode.

\begin{lemma}\label{lem Markov 2.1}
The triple $(X(\cdot),L(\cdot),R(\cdot))$ is a $(2m+2)$-dimensional
Markov process on every time interval $[\tau_i,\tau_{i+1})$, for
$i=0,1,\cdots,N-1$.
\end{lemma}
\preuve{}{ Denoting $\mathbf{1}=(1,1,\cdots,1)$ as the
$(m+1)$-dimensional row vector of 1's, and
$\mathbf{0}=(0,\cdots,0)$ as the $m$-dimensional row vector of
0's. Over the time interval $(\tau_i,\tau_{i+1})$, the triple $(X(\cdot),L(\cdot),R(\cdot))$ constitutes a
strong solution to the $(2m+2)$-dimensional SDE
\begin{equation}\label{sde 1.4}
\left\{ \begin{aligned}
dX(t)=&~\sigma(t,X(t))dW^0(t); \\
dL(t;\mu_j)=&~L(t;\mu_j)\frac{b(t,X(t);\mu_j)}{\sigma(t,X(t))}dW^0(t)\,
,\; \quad j=0,1,\cdots,m\,;\\
dR(t;\mu_j)=&~\frac{L(t;\mu_0)}{L(t;\mu_j)}\lambda e^{-\lambda
t}dt\,,\; \quad j=1,\cdots,m
                          \end{aligned} \right.
                          \end{equation}
driven by the standard $\prob^0$-Brownian motion $W^0(\cdot)$, with the initial value
\begin{equation}\label{sde 1.1.1}
(X(0),L(0),R(0))=(x_0,\mathbf{1},\mathbf{0})
\end{equation}
at time 0 and the initial value $(X(\tau_i),L(\tau_i),R(\tau_i))$ at
the time $\tau_i$. From Assumption \ref{assump sde 1.1} (i)(ii) and inequality (\ref{cond linear}), the coefficients of the
SDE (\ref{sde 1.4}) are bounded on compact subsets of $\re^{2m+2}$
and are locally Lipschitz. The SDE (\ref{sde 1.4}) has a pathwise
unique, strong solution. The well-posedness of the SDE (\ref{sde
1.4}) (equivalently, the well-posedness of the associated martingale
problem) implies the $\prob^0$-strong Markov property of
$(X(\cdot),L(\cdot),R(\cdot))$, with respect to the Borel
$\sigma-$algebra $\mathbb{F}$ (\cite{StrVar} Stroock and Varadhan
(1997)). But the filtration $\filtr$ generated by $X(\cdot)$ is
contained in $\mathbb{F}$, and the process
$(X(\cdot),L(\cdot),R(\cdot))$ is $\filtr$-adapted. Then
$(X(\cdot),L(\cdot),R(\cdot))$ has the strong Markov property under
the probability measure $\prob^0$ with respect to $\filtr$.}

\begin{proposition}\label{prop non explosion}
The solution $(X(\cdot),L(\cdot),R(\cdot))$ to the SDE (\ref{sde
1.4}) does not explode within the time horizon $[0,T]$.
\end{proposition}
\preuve{}{By the definition of the explosion time of an SDE with
locally Lipschitz coefficients (e.g. page 330 of \cite{KarShr88}
Karatzas and Shreve (1988)), this follows from Lemma \ref{lem UI
1.1}.}

\noindent Eventually, we are able to
reformulate in Theorem \ref{thm contr ref meas} the partially
observable impulse control problem (\ref{opt value 2.1.0}). Under
the reference probability measure $\prob^0$, it becomes a fully
observable impulse control problem of the $(2m+2)$-dimensional
$\filtr$-adapted state process $(X(\cdot),L(\cdot),R(\cdot))$ from
the SDE (\ref{sde 1.4}). To prove this theorem, we use equations
(\ref{opt value 2.1.0}), (\ref{opt value 2.1}) and Lemma \ref{lem
opt value 2.1}.

\begin{theorem}\label{thm contr ref meas}
The impulse control problem (\ref{opt value 2.1.0}) under the
physical measure $\prob$ is equivalent to an impulse control problem
under the reference probability measure $\prob^0$, by choosing an
optimal $(\tau^*,\zeta^*)=\{(\tau_i^*,\zeta_i^*)\}_{i=1}^N$  to
achieve the maximal expected reward
\begin{equation}
\label{opt value 2.3.1} V^0 := \hspace{-1mm}
\sup\limits_{(\tau,\zeta)\in\textbf{I}}\expect^0 \hspace{-1mm} \left[\int_0^T
\hspace{-2mm} \alpha\left(s,X(s),L(s), R(s)\right)ds+\hspace{-1mm}\sum\limits_{i=1}^N
\beta(\tau_i,X(\tau_i-),L(\tau_i-),R(\tau_i-),\zeta_i)\right].
\end{equation}
Furthermore, the two maximal expected rewards are related by $V =
\xi(x_0)+ V^0$.
\end{theorem}

Because the best expected values $V$ and $V_0$ are
different only up to a constant $\xi(x_0)$, the two suprema in
(\ref{opt value 2.1.0}) and (\ref{opt value 2.3.1}) are achieved by
the same set of optimal control $(\tau^*,\zeta^*)$, if any. The
impulse control problem (\ref{opt value 2.3.1}) is the one we shall
solve.

\subsection{Solution under the reference probability measure}\label{subsec
solution}

This subsection will solve the impulse control problem (\ref{opt
value 2.3.1}), by representing the optimal control
$(\tau^*,\zeta^*)$ in Proposition \ref{prop optimal contr} in terms
of the value function and the state process. The cornerstone of the
representation is the dynamic programming principle of Lemma
\ref{lem DPP 2.1}. To satisfy the technical condition of the Snell
envelope argument for Proposition \ref{prop optimal contr}, the
continuity of
the value functions is provided in Lemma \ref{lem cont}.\\
\\To save notations, some abbreviations are introduced first. We denote by $\open$ the range of the solution
$(X(\cdot),L(\cdot),R(\cdot))$, and its boundaries as
\begin{equation}
Q:=[0,T]\times \open \quad \text{ and } \quad \partial^* Q:= \{T\}\times \open.
\end{equation}
The state space $\open$ differs for different parameters
$b(\cdot,\cdot;u)$ and $\sigma(\cdot,\cdot)$. Without loss of
generality, the variational
inequalities associated with the impulse control problem shall be
studied over the largest possible domain, which is
\begin{equation}
\open=\re\times(0,\infty)^{m+1}\times [0,\infty)^{m}.
\end{equation}
For every $n=1,2,\cdots$, denote the bounded domain
\begin{equation}
\open_n:=(-n,n)\times\left(\frac{1}{n},n\right)^{m+1}\times\left[0,n\right)^m\subset \open\subset
\re^{2m+2}.
\end{equation}
The closure of $\open_n$, denoted as $\bar{\open}_n$, is strictly contained in $\open$. As $n\rightarrow \infty$, the  sets $\open_n$ increase to $\open$, hence the sets $Q_n:=[0,T]\times \open_n$ increase to $Q$.    We introduce the abbreviations
\begin{equation}
y=(x,l_0,l_1,\cdots,l_m,r_1,\cdots,r_m),
\end{equation}
\begin{equation}
b_Y(t,y)=\left(0,0,0,\cdots,0,\frac{l_0}{l_1}\lambda e^{-\lambda
t},\cdots,\frac{l_0}{l_m}\lambda e^{-\lambda t}\right),
\end{equation}
\begin{equation}
\sigma_Y(t,y)=\left(\sigma(t,x),l_0\frac{b(t,x;\mu_0)}{\sigma(t,x)},
l_1\frac{b(t,x;\mu_1)}{\sigma(t,x)},\cdots,l_m\frac{b(t,x;\mu_m)}{\sigma(t,x)}
,0,\cdots,0\right),
\end{equation}
and
\begin{equation}
\Gamma(y,z)=(x+\gamma(x,z),l,r),
\end{equation}
for all $(t,y)=(t,x,l,r)$ in $Q$. With this notation, the SDE
(\ref{sde 1.4}), which has a pathwise unique, strong solution
\begin{equation}
Y(\cdot)=(X(\cdot),L(\cdot;\mu_0),L(\cdot;\mu_1),
\cdots,L(\cdot;\mu_m),R(\cdot;\mu_1),\cdots,R(\cdot;\mu_m)),
\end{equation}
can be written in the vector form
\begin{eqnarray}\label{sde 1.3}
\left\{\begin{array}{c} \displaystyle
dY(t)= b_Y(t,Y(t))dt+\sigma_Y(t,Y(t))dW^0(t),\;\tau_i< t< \tau_{i+1};\vspace{1mm}\\
\displaystyle Y(\tau_i)= \Gamma(Y(\tau_i-),\zeta_i)\text{, for
}i=1,2,\cdots,N.
\end{array}\right.
\end{eqnarray}
Here the initial value is
\begin{equation}\label{sde 1.2.1}
Y(0)=(x_0,\mathbf{1},\mathbf{0}),
\end{equation}
and $W^0(\cdot)$ is a standard
$\prob^0$-Brownian motion.\\
\\In the abbreviated notation, the maximal expected reward in
equation (\ref{opt value 2.3.2}) can be written as
\begin{equation}
\label{opt value 2.3.2} V^0 =
\sup\limits_{(\tau,\zeta)\in\textbf{I}}\expect^0\left[\int_0^T
\alpha\left(s,Y(s)\right)ds+\sum\limits_{i=1}^N
\beta(\tau_i,Y(\tau_i-),\zeta_i)\right].
\end{equation}
The rest of this section will use the above abbreviated notations.

\begin{lemma}\label{lem DPP 2.1}{\bf Dynamic Programming Principle.}
For any $k\in \{1,2,\cdots,N\}$, and any
$0\leq t\leq T$, let $\textbf{I}_{t,k}$ be the set of admissible
interventions $\{(\tau_i,\zeta_i)\}_{i=N-k+1}^N $ such
that $\tau_{N-k+1}\geq t$. Suppose the current value of the state
process $Y(t)=y\in \open$. There exist
deterministic measurable functions $\, v_0,v_1,$
$\cdots,v_N:Q\rightarrow\re$, such that
\begin{equation}\label{opt value 2.5}
\begin{split}
v_k(t,y) =\hspace{2mm}\,\text{ess}\hspace{-12.2mm}
\sup\limits_{\{(\tau_i,\zeta_i)\}_{i=N-k+1}^N\in\textbf{I}_{t,k}}
\expect^0\left[\left.\int_t^T \alpha\left(s,Y(s)\right)ds
+\sum\limits_{i=N-k+1}^N\beta
(\tau_i,Y(\tau_i-),\zeta_i)\right|\filtr (t)\right],
\end{split}
\end{equation}
for $k=1,\cdots,N$, and
\begin{equation}\label{opt value 2.5.1}
\begin{split}
v_0(t,y) =\expect^0\left[\left.\int_t^T \alpha\left(s,Y(s)\right)ds
\right|\filtr (t)\right].
\end{split}
\end{equation}
The value functions $v_1,\cdots,v_N$ satisfy the dynamic programming
principle
$$
v_k(t,y) =\hspace{2mm}\,\text{ess}\hspace{-12.2mm}\sup\limits_{(\tau_{N-k+1},\zeta_{N-k+1})
\in\textbf{I}_{t,1}}\expect^0 \Bigg[ \int_t^{\tau_{N-k+1}} \alpha\left(s,Y(s)\right)ds
$$
\begin{equation}
\label{DPP 2.1}
+\beta
(\tau_{N-k+1},Y(\tau_{N-k+1}-),\zeta_{N-k+1}) + \,v_{k-1}
(\tau_{N-k+1},\Gamma(Y(\tau_{N-k+1}),\zeta_{N-k+1}))
 \Bigg|\filtr (t)\Bigg].
\end{equation}
\end{lemma}
\preuve{}{ The existence of the functions $v_0,v_1,\cdots,v_N$
comes from the Markovian structure of the state process $Y(\cdot)$,
by Lemma \ref{lem Markov 2.1}.

To prove the equation (\ref{DPP 2.1}), fix an arbitrary $k\in
\{1,2,\cdots,N\}$, an arbitrary $t\in [0,T]$ and arbitrary admissible
interventions $\{(\tau_i,\zeta_i)\}_{i=N-k+1}^N \in
\textbf{I}_{t,k}$, we denote
\begin{equation}
\begin{split}
A_k(t):=&\int_t^{\tau_{N-k+1}} \alpha\left(s,Y(s)\right)ds+\beta
(\tau_{N-k+1},Y(\tau_{N-k+1}-),\zeta_{N-k+1});\\
B_k:=& \int_{\tau_{N-k+1}}^T
\alpha\left(s,Y(s)\right)ds+\sum\limits_{i=N-k+2}^N\beta
(\tau_i,Y(\tau_i-),\zeta_i).
\end{split}
\end{equation}
Then
\begin{equation}
\begin{split}
&\expect^0\left[\left.\int_t^T \alpha\left(s,Y(s)\right)ds
+\sum\limits_{i=N-k+1}^N\beta
(\tau_i,Y(\tau_i-),\zeta_i)\right|\filtr (t)\right]\\
=\,&\,\expect^0\left[\left.A_k(t)+\expect^0\left[B_k\left|\filtr(\tau_{N-k+1})\right.\right]\right|\filtr
(t)\right].
\end{split}
\end{equation}
On one hand, taking supremum over
$(\tau_{N-k+1},\zeta_{N-k+1})\in\textbf{I}_{t,1}$ on both sides of
the inequality
\begin{equation}\label{}
\begin{split}
\,&\,\expect^0\left[\left.A_k(t)+\expect^0\left[B_k\left|\filtr(\tau_{N-k+1})\right.\right]\right|\filtr
(t)\right]\\
\leq \,&\,\expect^0\left[\left.A_k(t)+v_{k-1}
(\tau_{N-k+1},\Gamma(Y(\tau_{N-k+1}),\zeta_{N-k+1}))\right|\filtr
(t)\right]
\end{split}
\end{equation}
shows $v_k(t,y)$ less than or equal to the right hand side of
(\ref{DPP 2.1}). On the other hand, the inequality
\begin{equation}\label{}
v_k(t,y)\geq\expect^0\left[\left.A_k(t)+\expect^0\left[B_k\left|\filtr(\tau_{N-k+1})\right.\right]\right|\filtr
(t)\right]
\end{equation}
implies
\begin{equation}\label{}
v_k(t,y)\geq \expect^0\left[\left.A_k(t)+v_{k-1}
(\tau_{N-k+1},\Gamma(Y(\tau_{N-k+1}),\zeta_{N-k+1}))\right|\filtr
(t)\right]
\end{equation}
and thus $v_k(t,y)$ greater than or equal to the right hand side of
(\ref{DPP 2.1}).

  See \cite{FlemmingSoner93} Fleming \& Soner (1993),
\cite{Krylov80} Krylov (1980)
or \cite{Pham09} Pham (2009) for a more
detailed account for the dynamic programming principle.}

\begin{lemma}\label{lem cont} The value functions $v_0,v_1,\cdots,v_N$
defined in (\ref{opt value 2.5})
and (\ref{opt value 2.5.1}) are continuous in $(t,y)\in Q$. Over the
compact set $\bar{\open}_n$, they admit
  moduli of continuity $\omega_n:[0,\infty)\rightarrow[0,\infty)$,
  uniformly for all $0\leq t\leq T$,
meaning that
\begin{equation}\label{cond Lip 1.1}
|v_k(t,y^1)-v_k(t,y^2)| \leq  \omega_n(||y^1-y^2||) \,,\quad \text{
for all } \quad (t,y^1), (t,y^2)\in [0,T]\times \bar{\open}_n.
\end{equation}
\end{lemma}
\preuve{}{
By the continuity of solutions to SDEs (Theorem 5.2 in Chapter II on
page 229 of \cite{Kunita82} Kunita (1982)), and by the continuity of
the function $\gamma$ given in Assumption \ref{assump rwd 1.1}
(iii), the unique strong solution to the controlled SDE (\ref{sde
1.4}) is continuous in its initial value $(t,Y(t))=(t,y)\in Q$. We
shall also use the continuity of the functions $\alpha$ and $\beta$
in equations (\ref{cross var 1.1}) and (\ref{cross var 2.1}),
 Assumption \ref{assump rwd 1.1} (ii)(iii) and the uniform
 integrability Lemma \ref{lem UI 1.1}. Inductively applying the proof of Proposition 2.2 in
\cite{JaiLamLap89} Jaillet, Lamberton and Lapeyre (1990) to $v_k$,
for $k=0,1,\cdots,N$, we know that the value functions
$v_0,v_1,\cdots,v_N$ are continuous in $(t,y)\in Q$.

Restricted on the compact set $\bar{Q_n}=[0,T]\times \bar{\open}_n$,
the value functions $v_0,v_1,\cdots,v_N$ are uniformly continuous in
$(t,y) \in Q_n$, hence they admit moduli of continuity $\omega_n$ in the
space variable $y\in\bar{\open}_n$, for all $t\in[0,T]$. }

\noindent The collection of all continuous functions over the domain
$Q$, which admit the modulus of continuity $\omega_n$ for $y$ in the
compact set $\bar{\open}_n$ uniformly for all $0\leq t\leq T$, is
denoted as $\cont(Q;\omega_n)$. It is the very set of properties
described in Lemma \ref{lem cont}.\\
\\The optimal impulse controls are then obtained in terms of the
value functions $v_0,v_1,\cdots,v_N$, and of the triple
$(X(\cdot),L(\cdot),R(\cdot))=Y(\cdot)$. The triple
$(X(\cdot),L(\cdot),R(\cdot))=Y(\cdot)$ of processes in (\ref{sde
1.4}), which is adapted to the filtration $\filtr$ generated by the
observation $X(\cdot)$, can be viewed as a ``sufficient statistic"
for the optimization problem (\ref{opt value 2.1.0}). This
``sufficient statistic" that the decision maker needs to monitor
remains the same for all cumulative reward functions $h(\cdot)$, all
impulse control costs $c(\cdot)$ and all terminal reward functions
$\xi(\cdot)$ in (\ref{opt value 2.1.0}).
\begin{proposition}\label{prop optimal contr} (Iterative procedure for optimization)
For any measurable function $f: Q\rightarrow \re$, define a mapping $\textbf{M}$
by
\begin{equation}\label{mapping max 2.1}
\left(\textbf{M}f\right)(t,y)\,:=\,\sup\limits_{z\in\re}
\{f(t,\Gamma(y,z))+\beta(t,y,z)\}\,,~~~\text{ for all } \,(t,y)\in
Q.
\end{equation}
For every $k=1,2,\cdots,N$, iteratively define an $\filtr$-stopping
time
\begin{equation}\label{tau star}
\tau_k^*:=\inf \left\{\tau_{k-1}^*<t\leq T\left|
v_{N-k+1}\left(t,Y(t)\right)\leq \textbf{M}v_{N-k}\left(t,Y(t)\right)
\right. \right\},
\end{equation}
with the convention that $\tau_0^*=0$. Suppose the supremum
\begin{equation}\label{mapping max 2.2}
\sup\limits_{z\in\re}
\left\{v_{N-k+1}(t,\Gamma(y,z))+\beta(t,y,z)\right\}-v_{N-k}(t,y)
\end{equation}
can be attained by a real number $z_k(t,y)$, and define an $\filtr
(\tau_k^*-)$-measurable random variable
\begin{equation}
\zeta_k^*:=z_k\left(\tau_k^*,Y(\tau_k^*-)\right),
\end{equation}
for every $k=1,2,\cdots,N$. Then the suprema in (\ref{opt value
2.1.0}) and (\ref{opt value 2.1}) are attained by the set of impulse
controls $\left\{\tau_k^*,\zeta_k^*\right\}_{k=1}^N$ in $\textbf{I}$.
Furthermore, the maximal expected rewards $V^0=v_N(0,Y(0))$ and
$V=\xi(x_0)+v_N(0,Y(0))$.
\end{proposition}

\begin{remark}
When Bensoussan and Lions were originally formulating the impulse
control problem in the 1970's, their number of interventions
$N=\infty$. There is no fundamental difference whether $N$ is finite
or infinite, except that slightly different technical conditions on
the coefficients and the admissible control set are required to
derive properties like well-posedness, continuity and even
differentiability of the value function. It has been pointed out by
Bensoussan and Lions in Theorem 4 of \cite{BensLions75b} that the
value function of $N$ interventions converges to that of infinitely
many interventions, as $N\rightarrow\infty$. To extend results in
this paper to $N=\infty$ means modifying the technical assumptions.
\end{remark}

\section{Geometric Brownian motion with drift
uncertainty}\label{subsec GBM}

In this part, we discuss how to approximate the optimal stopping
time distribution $\prob (\tau^*_i \in (t, t+dt])$ thanks to a Monte
Carlo simulation, and use the results for constructing a trading
strategy. Thus, we start by introducing in Section \ref{sec41} the
example on which we implement the theoretical setting of Sections
\ref{subsec measure change} and \ref{subsec solution}. Then, we
present in Section \ref{sec42} a method based on Longstaff-Schwartz
algorithm to simulate the optimal stopping times $ \left\{
\tau^*_i\right\}_{i=1,...,N}$ family. In Section \ref{sec43}, we
give a simple static trading strategy that allows to test the
simulation results.

\subsection{Setting the problem \label{sec41}}
To illustrate the model (\ref{SDE impulse 2.1}), we discuss Geometric
Brownian Motion as a commonly seen simple example. The parameter $\theta(\cdot)$ is the
drift with the initial value $\mu_0$. The random variable $U$ has
the prior distribution
\begin{equation} \label{theta prior 1.2.1}
U=\left\{ \begin{aligned}
         \mu_1 &\text{, with probability }p_1;\\
         \mu_2 &\text{, with probability }p_2=1-p_1
                          \end{aligned} \right.
                          \end{equation}
and $\rho$ has an exponential $\lambda$ prior distribution as in
(\ref{rho prior 1.1}). The diffusion $X(\cdot)$ in (\ref{SDE impulse 2.1}) is the geometric Brownian
motion
\begin{eqnarray}\label{GBM 1.2.1}
\left\{ \begin{aligned} dX(t) & = X(t)\theta(t)dt+X(t)\sigma dW(t);\\
X(0)&= x_0.
\end{aligned} \right.
\end{eqnarray}
In this example, the volatility $\sigma$ is a deterministic positive
number. The parameter $\theta(\cdot)$ with the initial value $\mu_0$
is the percentage drift of the Geometric Brownian motion.

 Suppose $X(\cdot)$ is the price process of a certain stock, and
there is zero interest rate, no transaction cost and no price impact.
Observing the price evolution only, an optimal trading problem is
finding two stopping times $0\leq \tau_1^*\leq\tau_2^*\leq T$
in $\xst$, to achieve the supremum in
\begin{equation}\label{opt value 2.8}
\sup\limits_{\tau_1\text{ and }\tau_2\in \xst\text{,
}\tau_1\leq\tau_2} \expect\left[X(\tau_2)-X(\tau_1)\right].
\end{equation}
In terms of the money received, the value (\ref{opt value 2.8}) is the best possible average
profit from first buying, then selling, one share of this stock. Comparing
the SDEs (\ref{GBM 1.2.1}) and (\ref{SDE impulse 2.1}), and the goal
(\ref{opt value 2.1.0}) with (\ref{opt value 2.8}), we are trying to
solve the impulse control problem with $\gamma(\cdot,\cdot)=0$,
$h(\cdot)=0$ and $\xi(\cdot)=0$.  We should set $c(x,z)=zx$ for
$x\in\re^n$ and $z\in\re$, $\zeta_1=-1$, $\zeta_2=1$ and $N=2$. To
incorporate  transaction costs and   price impact, it only remains
to modify the functions $c(\cdot,\cdot)$ and $\gamma(\cdot,\cdot)$.

For the geometric Brownian motion example, we may compute to get the likelihood ratio
processes
\begin{equation}\label{lr 1.2.4}
L(t;u)=\exp \left\{\frac{1}{2}\left(u-\frac{u^2}{\sigma^2}\right)
t\right\} \left(\frac{X(t)}{x_0}\right)^{u/\sigma^2},\; \quad 0\leq t\leq T.
\end{equation}
The process $R(\cdot;u)$ defined in (\ref{lr 1.2}) can then be written as
\begin{eqnarray}\label{int 1.2.3}
\hspace{-7mm}R(t;u) \hspace{-0.5mm} = \hspace{-1.5mm} \int_0^t \hspace{-2mm} \lambda \exp
\left\{\hspace{-1mm}\left(\frac{1}{2}\left(\mu_0-u-\frac{\mu_0^2-u^2}{\sigma^2}\right)
-\lambda\right)s\right\}\hspace{-1mm}\left(\frac{X(s)}{x_0}\right)^{(\mu_0-u)/\sigma^2}
\hspace{-12mm}ds, \quad 0\leq t\leq T,
\end{eqnarray}
for $u\in\Theta$. The dimensionality of the variational inequality can be
reduced, by using this alternative expression of $L(t;u)$ in terms of $X(t)$.
Substituting $\displaystyle l_j=\exp
\left\{\frac{1}{2}\left(\mu_j-\frac{\mu_j^2}{\sigma^2}\right)
t\right\}\left(\frac{x}{x_0}\right)^{\mu_j/\sigma^2}$, $\gamma(\cdot,\cdot)=0$, $h(\cdot)=0$,
$\xi(\cdot)=0$ and $c(x,z)=zx$ in two the functions $\alpha$ and $\beta$
defined in (\ref{cross var 1.1}) and (\ref{cross var 2.1}), we get
\begin{equation}\label{cross var 1.2.3}
\begin{split}
\alpha(t,x,l,r) = 0,
\end{split}
\end{equation}
\begin{eqnarray}\label{cross var 1.2.3}
\begin{array}{ccc}
\beta(t,x,l,r,z) & = & z\left( \begin{array}{c} \displaystyle \sum\limits_{j=1}^2
p_j\exp \left\{\frac{1}{2}\left(\mu_j-\frac{\mu_j^2}{\sigma^2}\right)
t\right\}\frac{x^{1+\mu_j/\sigma^2}}{x_0^{\mu_j/\sigma^2}} \hspace{1mm} r_j \vspace{1mm} \\ \displaystyle + \exp
\left\{\frac{1}{2}\left(\mu_0-\frac{\mu_0^2}{\sigma^2}\right)
t-\lambda t\right\}\frac{x^{1+\mu_0/\sigma^2}}{x_{0}^{\mu_0/\sigma^2}} \end{array} \right)
\vspace{3mm} \\ & =: & \hspace{-5.7cm} \bar{\beta}(t,x,r,z).
\end{array}
\end{eqnarray}

Under the measure $\prob^0$, the supremum in
(\ref{opt value 2.8}) becomes
\begin{equation}\label{opt value 2.9}
\begin{split}
\sup\limits_{\tau_1\text{ and }\tau_2\in \xst\text{,
}\tau_1\leq\tau_2} \expect&\left[X(\tau_2)-X(\tau_1)\right]\\
=\sup\limits_{\tau_1\text{ and }\tau_2\in \xst\text{,
}\tau_1\leq\tau_2}
\expect^0&\left[\bar{\beta}(\tau_2,X(\tau_2),R(\tau_2),1)+
\bar{\beta}(\tau_1,X(\tau_1),R(\tau_1),-1)\right].
\end{split}
\end{equation}
There exist deterministic measurable functions $\bar{v}_1$ and
$\bar{v}_2:[0,T]\times(0,\infty)\times
[0,\infty)^2$, such that
\begin{equation}
\begin{split}
\bar{v}_1(t,X(t),R(t))=&\sup\limits_{\tau_2\in\st_t}
\expect^0\left[\bar{\beta}(\tau_2,X(\tau_2),R(\tau_2),1)\left|\filtr (t)\right.\right];\\
\bar{v}_2(t,X(t),R(t))=&\sup\limits_{\tau_1\in\st_t}
\expect^0\left[\bar{\beta}(\tau_1,X(\tau_1),R(\tau_1),-1) +
\bar{v}_1(\tau_1,X(\tau_1),R(\tau_1))\left|\filtr (t)\right.\right].
\end{split}
\label{v1v2}
\end{equation}
The optimal value of the round-way transaction is
\begin{equation}\label{opt value 2.9}
\sup\limits_{(\tau_1,\tau_2)\in \xst^2 \atop \tau_1\leq\tau_2}
\expect\left[X(\tau_2)-X(\tau_1)\right]=\bar{v}_2(0,X(0),\mathbf{0}).
\end{equation}

\subsection{Longstaff-Schwartz algorithm for multiple optimal stopping times\label{sec42}}

Proposed in \cite{long} for pricing American options, the Longstaff-Schwartz procedure was rigourously
formulated by Cl\'ement, Lamberton and Protter (2002) in \cite{lamb} using stopping times instead of the value function
for the dynamic programming algorithm. This method also involves a regression approximation of the conditional
expectation and the discretization of the interval on which the stopping times take their values. The convergence
due to each approximation step is studied in \cite{lamb} and we aim at reusing this algorithm for our multiple
optimal stopping problem.

Applied also by Tsitsiklis and Van Roy (2001) in \cite{tsi}, the
regression approximation of the conditional expectation was extended
to BSDEs (Backward Stochastic Differential Equations) by Gobet,
Lemor and Warin (2005) in \cite{gobet}. The regression approximation generally uses one regression
vector for the whole set of trajectories and it is then considered
as a global method. Thus, other authors use more
local approximations of the conditional expectation,
based on either Malliavin calculus as in \cite{Lok2,bouch} or
quantization method as in \cite{pages}. To keep the presentation of
our algorithm simple, we apply a monomial regression method.
However, one has to keep in mind that in some problems, especially
when the dimension becomes high (more than two assets), a good
approximation of the conditional expectation is a key ingredient.

To proceed, we first need to approach stopping times in $\xst $ with stopping times taking values in
the finite set $0=t_{0}<t_{1}<...<t_{n}=T$. Then, the computation of
\eqref{v1v2} can be reduced to the implementation of two dynamic programming algorithms that
we express in terms of the optimal stopping times $\tau_1^k$ and $\tau_2^k$,
for each path, as follows
\begin{eqnarray}
\begin{array}{c}
\displaystyle \tau_2^n=T,\vspace{2mm}\\
\displaystyle \text{for } k \in \{ n-1,...,1 \},\quad \tau_2^k = t_k 1_{A_{1k}} + \tau_2^{k+1} 1_{A^{c}_{1k}}, \label{tauTwo}
\end{array}
\end{eqnarray}
\begin{eqnarray}
\begin{array}{c}
\displaystyle \tau_1^n =  T,\vspace{2mm}\\
\displaystyle \text{for } k \in \{ n-1,...,1 \},\quad \tau_1^k = t_k 1_{A_{2k}} + \tau_2^k \wedge \tau_1^{k+1} 1_{A^{c}_{2k}} \label{tauOne}
\end{array}
\end{eqnarray}
and, denoting $\expect^0_{t_k}$ the conditional expectation knowing $\filtr (t_{k})$, the sets $A_{1k}$ and $A_{2k}$ are given by
\begin{eqnarray*}
\hspace{-3mm}\begin{array}{c}
A_{1k} \hspace{-1mm} = \hspace{-1mm} \left\{ \bar{\beta}(t_k,X(t_k),R(t_k),1) >
\expect^0_{t_k} \hspace{-1mm} \left[\bar{v}_1(t_{k+1},X(t_{k+1}),R(t_{k+1})) \right] \right\};
\vspace{2mm}\\ A_{2k} \hspace{-1mm} = \hspace{-1mm} \left\{ \begin{array}{c} \bar{\beta}(t_k,X(t_k),R(t_k),-1) \vspace{1mm} \\ +
\bar{v}_1(t_{k},X(t_{k}),R(t_{k})) \end{array} > \expect^0_{t_k} \hspace{-1mm} \left[\bar{v}_2(t_{k+1},X(t_{k+1}),R(t_{k+1})) \right] \right\}.
\end{array}
\label{A1A2}
\end{eqnarray*}

While the simulation of $X$ under $\prob^0$ is straightforward, the simulation
of $R$ is performed thanks to a trapezoidal approximation of the
time integral applied in \cite{lap} for Asian options. It remains
then to compute $\expect^0_{t_k} \hspace{-1mm}
\left[\bar{v}_i(t_{k+1},X(t_{k+1}),R(t_{k+1})) \right]$ for $i=1,2$.
Employing the Markov property established in Lemma \ref{lem Markov
2.1}, a conditional expectation according to $\filtr (t_{k})$ can be
replaced by a conditional expectation according to
$Y(t_{k})=(X(t_{k}),R(t_{k}))$. In our application, this latter
quantity will be approximated by a regression on the monomial family
$g (Y(t_{k}))= (1,X(t_{k}),R(t_{k},\mu_1),R(t_{k},\mu_2))$.
Formally, for the auxiliary functions $f_i (\cdot ) = \bar{v}_i(t_{k},\cdot,\cdot)$,
$i=1,2$
\begin{eqnarray}
\expect^0 \left( f_i(Y(t_{k+1}))| Y(t_{k}) \right)
\approx A^i.g (Y(t_{k})).
\end{eqnarray}
The vector $A^i$ minimizes the quadratic error
\begin{eqnarray}
\left| \left| f_i(Y(t_{k+1})) - A^i.g (Y(t_{k})) \right| \right|_{L^{2}} \label{sqrtCH1}
\end{eqnarray}
and thus equal to
\begin{eqnarray}
A^i = \Psi^{-1} \expect^0 \left( f_i(Y(t_{k+1})) g (Y(t_{k})) \right), \label{valCH1}
\end{eqnarray}
where the matrix $\Psi = \expect^0 \left( g (Y(t_{k})) g^{t} (Y(t_{k})) \right)$ and $^{t}$ is the transpose
operator. Consequently, at each time step, the matrix inversion (\ref{valCH1}) can be implemented by the Singular
Value Decomposition (SVD) explained in \cite{NumRes} and the expectations are approximated by an arithmetic average
\begin{eqnarray*}
\expect^0 \left( f_i(Y(t_{k+1})) g (Y(t_{k})) \right) \approx \frac{1}{M}\sum_{l =1}^{M}
f_i(Y^l(t_{k+1})) g (Y^l(t_{k})),
\end{eqnarray*}
\vspace{-5mm}
\begin{eqnarray*}
\expect^0 \left( g (Y_{t_{k}}) g^{t} (Y_{t_{k}}) \right) \approx \frac{1}{M}\sum_{l =1}^{M} g (Y^l(t_{k}))
g^t (Y^l(t_{k})).
\end{eqnarray*}
and $M$ is the number of simulated trajectories. Then, using $\tau^1_1$ known from \eqref{tauOne}
\begin{equation}\label{appr value 2.9}
\bar{v}_2(0,X(0),0) \approx
\max \left(\expect^0\left[\bar{v}_2(\tau^1_1,X(\tau^1_1),R(\tau^1_1))\right],0 \right)
\end{equation}
is the approximation of \eqref{opt value 2.9}. Also, for $i=1,2$ we make the approximation
\begin{equation}\label{appr tau0}
\prob^0 (\tau^*_i \in (t_{k}, t_{k+1}]) \approx
\prob^0 (\tau^1_i = t_{k+1}).
\end{equation}
Finally, we should point out that the proposed procedure can be generalized for more than
two optimal stopping times. The convergence of the overall algorithm can be established in the
same way as it is presented in \cite{lamb} for one optimal stopping time. Besides, the reader
should notice that we only considered the deterministic case $\zeta_1=-1$, $\zeta_2=1$ and $N=2$.
In fact, one can propose a randomized version of the algorithm proposed above that includes an
optimization over $\zeta $, however our purpose here is only to give an illustration of a simple case.
As a future work, we will study the convergence of a more general method for impulse control
based on the multiple optimal stopping times algorithm presented above.

\subsection{A trading strategy based on $\prob (\tau^*_i \in (t, t+dt])$\label{sec43}}
To present this trading strategy, we need first to change the probability measure and go back
to $\prob $ thanks to \eqref{change measure 1.1} and the simulation of $Z(t)$ using
\eqref{RN deriv 1.2}. For $i=1,2$, we obtain then the approximation
\begin{equation}\label{appr tau}
\prob (\tau^*_i \in (t_{k}, t_{k+1}]) \approx
\expect^0 \left( Z(t_{k+1}) 1_{\tau^1_i = t_{k+1}} \right).
\end{equation}

Now, let us assume that we can buy and sell not only one stock but a
bigger volume $q \geq 1$ of stocks. Consequently, one can use the
approximation
\begin{equation}\label{appr val}
q \hspace{-3mm}\sup\limits_{(\tau_1,\tau_2)\in \xst^2 \atop \tau_1\leq\tau_2} \hspace{-3mm}
\expect\left[X(\tau_2)-X(\tau_1)\right] \approx q \max
\left(\expect^0\left[\bar{v}_2(\tau^1_1,X(\tau^1_1),R(\tau^1_1))\right],0 \right).
\end{equation}
Using the value obtained in \eqref{appr val}, we decide at a first stage if it is interesting
to trade or not. Indeed, if this value is not big ``enough'' then it is not worthwhile taking
the trading risks of losing money. In this one dimensional example, one should invest on $X$ only if
it is drifting more positively than negatively. Besides, if we are satisfied by the expected profits,
we can establish the following static trading strategy on the spot prices $X(t_{k+1})$ for $k= 0,..., n-1$:
The money received is
\begin{eqnarray}\label{MoneyR}
\hspace{-8mm} M^{r*}_{k+1}= M^{r*}_{k} + q\left[ \prob (\tau^*_2 \in (t_{k}, t_{k+1}]) -
\prob (\tau^*_1 \in (t_{k}, t_{k+1}])\right] X(t_{k+1}),\ M^{r*}_{0} = 0.
\end{eqnarray}
$M^{r*}$ can take negative values which mean that we are buying stocks. We are going to compare
$M^{r*}_{n}$ to $M^{r}_n$ where $M^{r}_{k+1}$ is defined by
\begin{equation}\label{MoneyRS}
M^{r}_{k+1}=M^{r}_{k} + q\left[ P^2_{k+1} - P^1_{k+1}\right] X(t_{k+1}),\ M^{r}_{0} = 0.
\end{equation}
and the quantities $P^1_{k+1}$, $P^2_{k+1}$ are simulated thanks to the following increasing
induction on $k= 0,..., n-1$
\begin{eqnarray}\label{MoneyRSS}
\begin{array}{c}
P^1_{k+1} \sim U(0,1-S^1_{k+1}),\quad
P^2_{k+1} \sim U(0,S^2_{k+1}),\quad P^2_{n} = S^2_{n}, \vspace{2mm}\\
S^1_{k+1} = S^2_{k}-P^2_{k} ,\quad
S^2_{k+1} = S^1_{k+1}+P^1_{k+1},\quad S^1_{1} = 0
\end{array}
\end{eqnarray}
where $U(0,x)$ is the uniform law on $[0,x]$.

Consequently, we are going to compare the money earned from our
static trading strategy \eqref{appr tau}\eqref{MoneyR} to some
$M_a^s$ arbitrary strategies specified by
\eqref{MoneyRS}\eqref{MoneyRSS}. We process this comparison on a
large number $M^{\text{new}}$ (given later) of newly
simulated trajectories of $X$ under the probability
$\prob$. In the following, we denote respectively by
$\widetilde{M^{r*}_{n}}$ and $\widetilde{M^{r}_n}$ the average value
of $M^{r*}_{n}$ and $M^{r}_n$ on the $M^{\text{new}}$ newly
simulated trajectories of $X$. Also we denote by
$\overline{M^{r*}_{n}}$ and $\overline{M^{r}_n}$ the maximum value
of respectively $M^{r*}_{n}$ and $M^{r}_n$ on the $M^{\text{new}}$
newly simulated trajectories of $X$.

\begin{figure}[t]
 \begin{minipage}[b]{.48\linewidth}
  \epsfig{figure=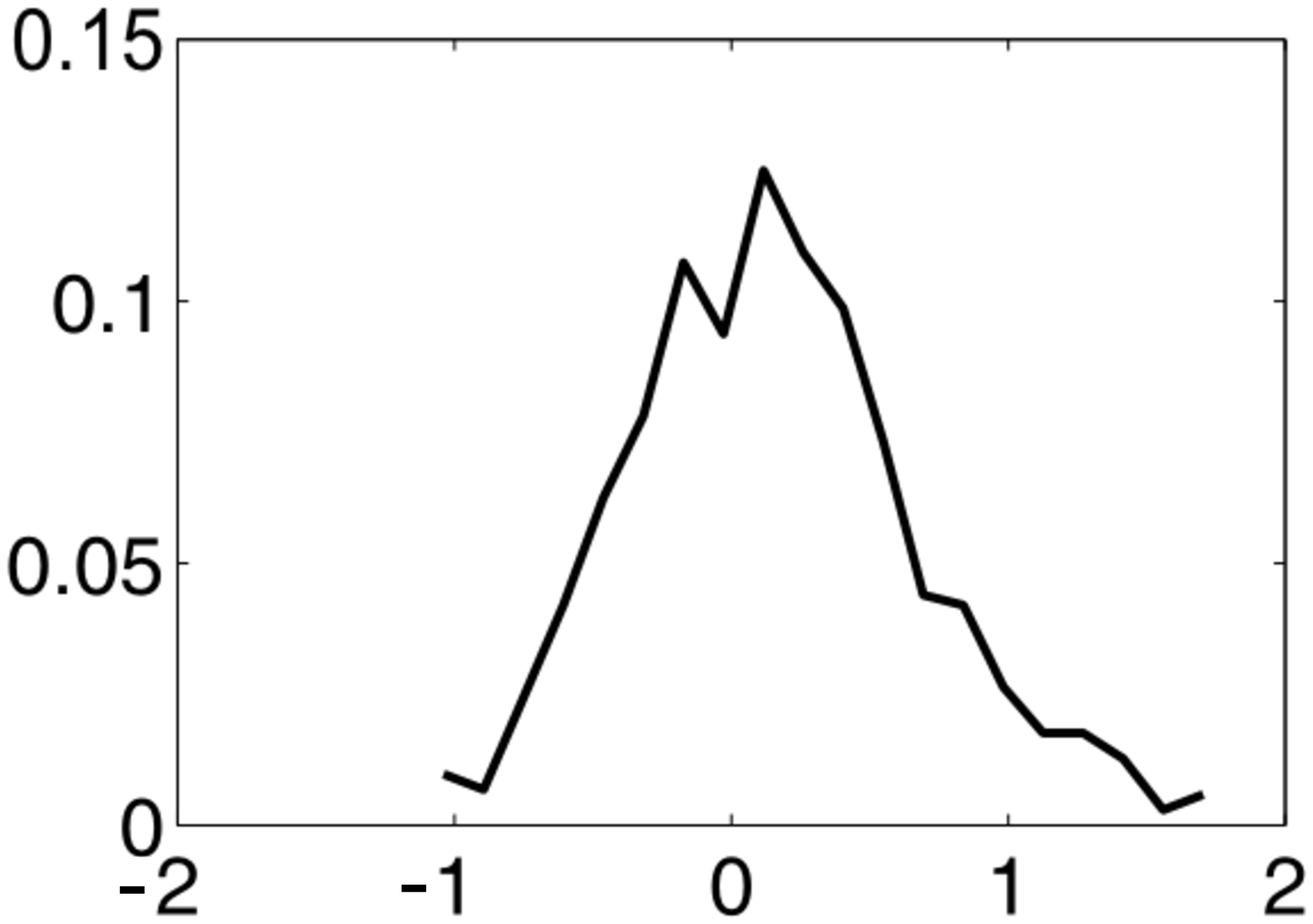,width=1.0\linewidth}
  \caption{\label{fig1} The histogram of $\widetilde{M^{r}_n}$ associated to
  $M_a^s = 2^{10}$ arbitrary strategies. The static optimal strategy provides $\widetilde{M^{r*}_n} = 6.17$.}
 \end{minipage} \hfill
 \begin{minipage}[b]{.48\linewidth}
  \epsfig{figure=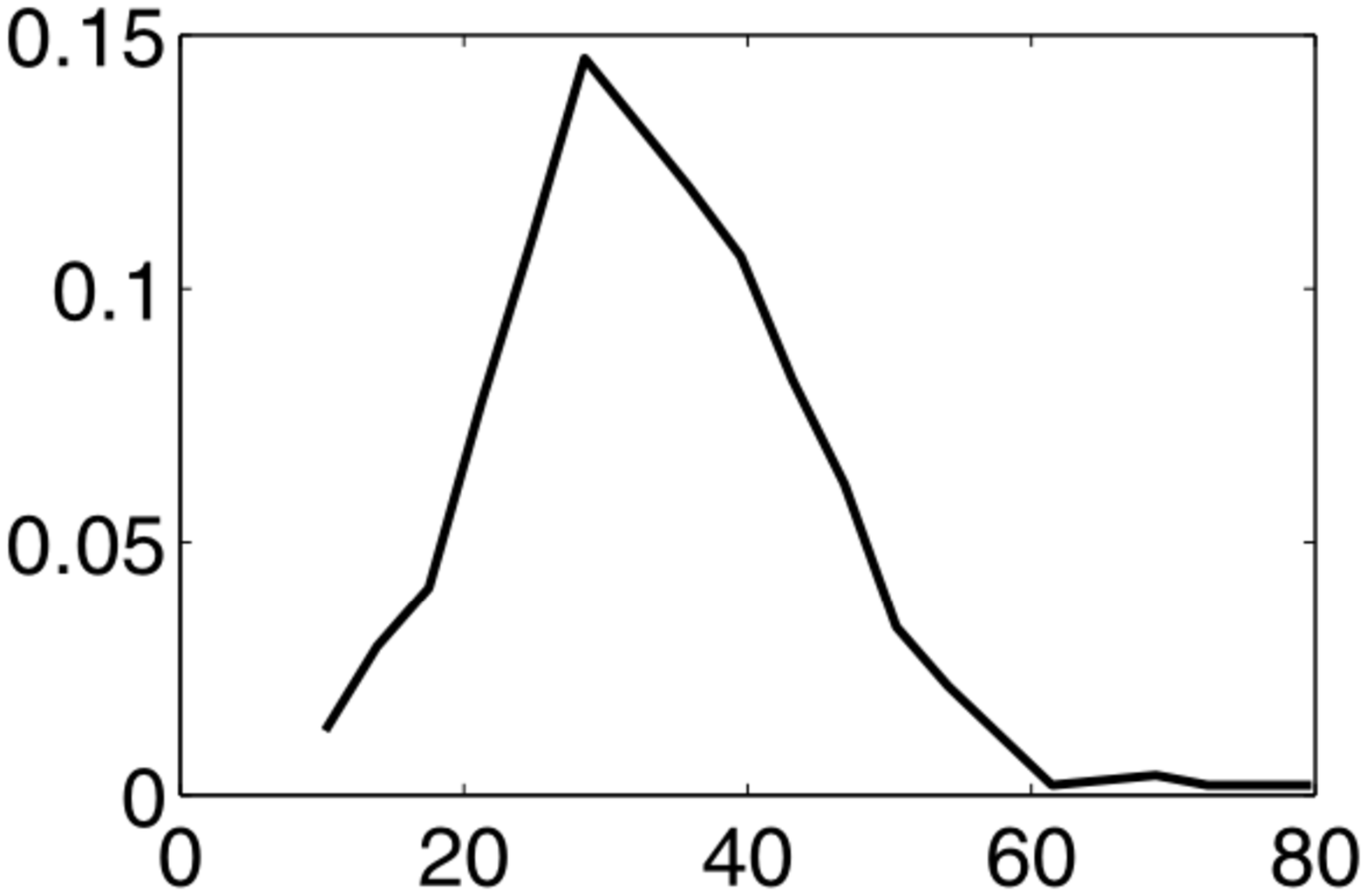,width=1.05\linewidth}
  \caption{\label{fig2} The histogram of $\overline{M^{r}_n}$ associated to
  $M_a^s = 2^{10}$ arbitrary strategies. The static optimal strategy provides $\overline{M^{r*}_n} = 51.98$.}
 \end{minipage}
 \end{figure}

\begin{figure}
 \epsfxsize=7cm
\centerline{\epsffile{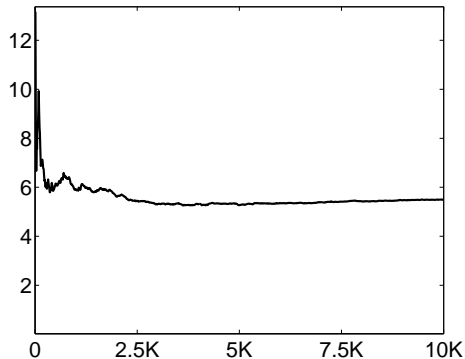}}
\caption{ \label{fig3} The evolution of $\widetilde{M^{r*}_n}$ according to the number of simulated trajectories $M^{\text{new}}$ (in kilo)}
\end{figure}

Although we tested our algorithm for a large number of model parameters, we present here the
results associated to only one choice of values. We refer the reader to the first author web page
to download the C++ code of the algorithm in order to test it with other parameters values.
Figures \ref{fig1}, \ref{fig2} and \ref{fig3} involve the following choice:
Number of simulated trajectories for Longstaff-Schwartz algorithm $M = 2^{16}$, number of
time steps $n= 10$, $T=1$, $\mu_0  = \mu_1  = 0.1$, $\mu_2  = -0.1$
$p_1= 0.5 $, $\sigma = 0.2 $, $ \lambda = 1 $, $ x_0 = 1$ and $q=100$.

Using $M^{\text{new}}=2^{10}$ new scenarios of the evolution of $X$,
we compare in figures \ref{fig1} and \ref{fig2} the average profit
as well as the maximum profit generated by the $M_a^s = 2^{10}$
arbitrary strategies to the ones generated by the static optimal
strategy. In Figure \ref{fig1}, the optimal strategy outperforms all
the arbitrary strategies which confirms the effectiveness of the
method implemented in Section \ref{sec42}. Moreover, even the
maximum profit provided by the optimal strategy is among the best
according to Figure \ref{fig2}. In Figure \ref{fig3}, we show the
stability of the optimal strategy to reach the average value of
profits even for small numbers of scenarios $M^{\text{new}} \sim
200$.

To conclude this section, although one can establish more elaborate trading strategies
using the approximated
values of $\prob (\tau^*_i \in (t, t+dt])$, the one that we provided in this section
allowed us
to show the efficiency of the Longstaff-Schwartz algorithm for multiple optimal stopping
times.

\section{Discussions\label{discu}}

In Section \ref{sec51}, we discuss how the change of measure method can be extended to the multidimensional
case. In Section \ref{sec52}, we briefly present the posterior probabilities method and
explain what make it difficult to apply to the multidimensional partially-observed control problems.

\subsection{The measure change method for multidimensional state processes\label{sec51}}
The measure change method proposed in Section \ref{sec general
theory 2.2} can be extended to the case when the diffusion in
Section \ref{subsec 2.1} is multidimensional, mostly by
replacing the notations for scalars to those for matrices. This
subsection will give the multidimensional version of the formulae
whose modifications are not very straightforward.

Suppose the diffusion $X(\cdot)$ in equation (\ref{sde 1.2}) becomes a
$(d\times 1)$-dimensional process driven by a $(d\times
1)$-dimensional Brownian motion $W^0(\cdot)$ with independent
components. Correspondingly, the coefficients should be modified
according to the dimensionality.

For any possible values $u\in\Theta$, the drift
$b(\cdot,\cdot\,;u):[0,T]\times \re^d \rightarrow \re^d$ is a
mapping valued in $\re^d$ and the volatility
$\sigma(\cdot,\cdot):[0,T]\times \re^d \rightarrow \re^{d\times d}$
is a $(d\times d)$-matrix-valued mapping. Let $||\cdot||$ denote the
Euclidean norm. Assumption \ref{assump sde 1.1} is replaced by
Assumption \ref{assump sde 5.1}.
\begin{assumption}
\label{assump sde 5.1}
There exists a constant $C>0$, such that\\
(i) for all $(t,x^1)$, $(t,x^2)\in [0,T]\times \re^d$, and for all
$u \in \Theta$, we have
\begin{equation}
\begin{split}
&||b(t,x^1;u)-b(t,x^2;u)||+||\sigma(t,x^1)-\sigma(t,x^2)||\\
&+\left|\left|\frac{b(t,x^1;u)}{\sigma(t,x^1)}-\frac{b(t,x^2;u)}{\sigma(t,x^2)}\right|\right|
\leq  C||x^1-x^2||\,;
\end{split}
\end{equation}
(ii) for all $(t,x)\in [0,T]\times \re^d$ and  all $u \in \Theta$,
the matrix $\sigma(t,x)$ is invertible and
\begin{equation}\label{cond bdd}
\left|\left|\frac{b(t,x;u)}{\sigma(t,x)}\right|\right|\leq C\,.
\end{equation}
\end{assumption}

\noindent The reward functions $\xi$ and $h:\re^d\rightarrow \re$,
as well as the intervention impact $\gamma$ and the reward from
intervention $c:\re^d\times\re\rightarrow\re$, satisfy Assumption
\ref{assump rwd 5.1} instead of Assumption \ref{assump rwd 1.1}.
\begin{assumption}
\label{assump rwd 5.1} (i) The function $\xi(\cdot)$ is twice
continuously differentiable, with first and second order derivatives
denoted as $\D{\xi}{x_i}(\cdot)$ and
$\CDD{\xi}{x_i}{x_j}(\cdot)$, for $i,j=1,\cdots,d$.\\
(ii) The functions $h(\cdot)$, $\xi(\cdot)$, $\D{\xi}{x_i}(\cdot)$
and $\CDD{\xi}{x_i}{x_j}(\cdot)$ are locally Lipschitz and have polynomial
growth, for $i,j=1,\cdots,d$.\\
(iii) The function $\gamma(x,z)$ is bounded for all $x\in\re^d$ and
$z\in\re$,
 and the function $c(x,z)$ has
polynomial growth rate in $x\in\re^d$ uniformly for all $z\in\re$.
Both functions $\gamma(x,z)$ and $c(x,z)$ are continuous in $x$, for
any arbitrarily fixed $z\in\re$.
\end{assumption}

\noindent To change between the physical measure $\prob$ and the
reference probability measure $\prob^0$, the likelihood ratio
process (\ref{lr 2.1}) is replaced by
\begin{equation}\label{lr 5.1}
L(t;u)=\exp \left\{ \begin{aligned} &\int_0^t
\left(\sigma^{-1}(s,X(s))b(s,X(s);u)\right)^{t}dW^0(s)\\
&-\frac{1}{2}\int_0^t \left|\left| \sigma^{-1}(s,X(s))b(s,X(s);u)
\right|\right|^2 ds \end{aligned}\right\},\;0\leq t\leq T\,.
\end{equation}
where $^{t}$ is the transpose operator. Then, by the same derivation
as in Section \ref{subsec measure change}, we arrive at the impulse
control problem (\ref{opt value 2.3.1}) under the reference
probability measure $\prob^0$. Let $\nabla$ denote the gradient
operator of a function. Instead of equations (\ref{cross var 1.1})
and (\ref{cross var 2.1}), the reward functions $\alpha$ and $\beta$
in (\ref{opt value 2.3.1}) are defined as
\begin{equation}\label{cross var 5.1}
\begin{split}
\alpha(t,x,l,r) :=&\left(\sum\limits_{j=1}^m p_jl_jr_j+e^{-\lambda
t}l_0\right)\left(h(x)+ \frac{1}{2}\sum\limits_{i,j=1}^d
\left(\sigma\sigma^{t}\right)_{i,j}(t,x)\CDD{\xi}{x_i}{x_j}(x)\right)
\\
&+\left((\nabla \xi)(x)\sum\limits_{j=1}^m
p_jl_jr_jb(t,x;\mu_j)+e^{-\lambda t}l_0b(t,x;\mu_0)\right),
\end{split}
\end{equation}
and
\begin{equation}\label{cross var 5.1}
\beta(t,x,l,r,z) := \left(\sum\limits_{j=1}^m p_jl_jr_j+e^{-\lambda
t}l_0\right)\left(\begin{aligned}&\xi(x+\gamma(x,z))-\xi(x)\\
&+(\nabla \xi)(x)\gamma(x,z)+c(x,z)
\end{aligned}\right).
\end{equation}
The solution to the impulse control problem (\ref{opt value 2.3.1})
with the multidimensional $X(\cdot)$ process will follow exactly the
same steps as in Section \ref{subsec solution}.

\subsection{Partial observation control via posterior probabilities\label{sec52}}

The traditional method to reduce a partially-observed control problem
to one with full observation is to augment the state process $X(\cdot)$ by
the posterior probability processes
\begin{equation}
\Pi_i(t):=\prob(\theta(t)=\mu_i|\filtr(t))\text{, for }
i=0,1,\cdots,m.
\end{equation}
This method is presented in Section 2.4.6 of \cite{Pham05}
Pham(2005) for a survey of the control problem and in Chapter 9 of
\cite{LiptserShiryaev2001} Lipster and Shiryaev (2001) for the
derivation of the posterior expectation and probabilities. In this
subsection, we shall first outline how to solve our problem in
dimension one by the posterior probability
method, modulus technical assumptions, and then briefly explore the relation between the two methods.\\
\\Define a function $\bar{b}: [0,T]\times \re \times
[0,1]^{m+1}\rightarrow \re$, $(t,x,\pi)\mapsto \bar{b}(t,x,\pi)$, by
$\bar{b}(t,x,\pi):=\sum\limits_{i=0}^m \pi_i b(t,x;\mu_i)$. Then the
uncertain drift projected onto the observation filtration is
\begin{equation}
\expect\left[ \left. b(t,X(t);\theta (t))\right| \filtr(t)\right]
=\bar{b}(t,X(t),\Pi(t)).
\end{equation}
Let $\bar{W}$ be the innovation Brownian motion. Given an
arbitrary admissible impulse control $(\tau,\zeta)\in\textbf{I}$, the
augmented state process $(X(\cdot),\Pi(\cdot))$ is a
$(m+2)$-dimensional Markov process on every time interval
$[\tau_k,\tau_{k+1})$, for $k=0,1,\cdots,N-1$, because it is the
unique strong solution to the controlled SDE
\begin{equation}\label{Newsde}
\left\{\begin{aligned} & X(t)=x_0+\int_0^t
\bar{b}(s,X(s),\Pi(s))ds+\int_0^t \sigma(s,X(s))d\bar{W}(s) +
\sum\limits_{\tau_i\leq
t}\gamma(X(\tau_i-),\zeta_i);\\
& \Pi_0(t)  = 1 - \lambda\int_0^t \Pi_0(s) ds +
\int_0^t \frac{b(s,X(s);\mu_0)-\bar{b}(s,X(s),\Pi(s))}{\sigma(s,X(s))}\Pi_0(s)d\bar{W}(s);\\
& \Pi_i(t)  = p_i\lambda \int_0^t \Pi_i(s) ds +
\int_0^t \frac{b(s,X(s);\mu_i)-\bar{b}(s,X(s),\Pi(s))}{\sigma(s,X(s))}\Pi_i(s)d\bar{W}(s),\\
&~~~~~~~~~~~~~~~~~~~~~~~~~~~~~~~~~~~~~~~~~~~~~~~~~~~~~~~~~~~
i=1,\cdots,m, 0\leq t\leq T.
\end{aligned}\right.
\end{equation}
One can then use the state process $(X(\cdot),\Pi(\cdot))$ to solve
the impulse control problem (\ref{opt value 2.1.0}) under the
physical measure, as a problem of full observation. The optimal impulse
controls are represented through the routine dynamic programming arguments
in terms of the value functions and the state process.\\
\\Let us proceed to demonstrate that the posterior probability method and the
measure change method are theoretically equivalent. Comparing with
those in Lemma \ref{lem DPP 2.1}, there exist value functions $\,
u_0,u_1,$ $\cdots,u_N:[0,T]\times \re \times [0,1]^{m+1}\rightarrow
\re$ such that
\begin{equation}
\begin{split}
u_k(t,x,\pi) =\hspace{2mm}\,\text{ess}\hspace{-12.2mm}
\sup\limits_{\{(\tau_i,\zeta_i)\}_{i=N-k+1}^N\in\textbf{I}_{t,k}}
\expect\bigg[&\int_t^T h\left(s,X(s)\right)ds + \xi(X(T))
\\&+\sum\limits_{i=N-k+1}^N c (X(\tau_i-),\zeta_i)\bigg|\filtr
(t)\bigg],
\end{split}
\end{equation}
for $k=1,\cdots,N$, and
\begin{equation}
\begin{split}
u_0(t,x,\pi) =\expect\left[\left.\int_t^T h\left(s,X(s)\right)ds +
\xi(X(T)) \right|\filtr (t)\right].
\end{split}
\end{equation}
By the same reasoning that derives Theorem \ref{thm contr ref meas},
the two sets of value functions respectively from the posterior
probability method and the measure change method are related by the
equations
\begin{equation}\label{post prob relation 5.1}
u_k(t,x,\pi)=\xi(x) + v_k(t,x,l,r)\text{, }k=0,1,\cdots,N,
\end{equation}
for all $x\in\re$, $\pi\in[0,1]^{m+1}$, $l\in(0,\infty)^{m+1}$ and
$r\in[0,\infty)^m$. In the case where the conditions in the Implicit
Mapping Theorem are satisfied, there exists an implicit mapping
$\bar{\pi}: Q\rightarrow [0,1]^{m+1}$, $(t,x,l,r)\mapsto
\bar{\pi}(t,x,l,r)$, such that
\begin{equation}\label{post prob relation 5.2}
u_k(t,x,\bar{\pi}(t,x,l,r))=\xi(x) + v_k(t,x,l,r)\text{,
}k=0,1,\cdots,N,
\end{equation}
for all $x\in\re$, $l\in(0,\infty)^{m+1}$ and $r\in[0,\infty)^m$.
The expression (\ref{post prob relation 5.2}) suggests that, when
applicable, the value functions before and after the change of
measure are different up to a change of variable.

Despite of the above equivalence, the measure change method has an
advantage in several dimensions when it comes to the numerical
implementation of Monte Carlo. Indeed, even when $X$ is
one-dimensional, one can easily remark that the Monte Carlo
simulation of \eqref{sde 1.4} is easier to perform and study than
the simulation of \eqref{Newsde}. With the latter SDEs system, one
has to propose an efficient discretization scheme and prove its
convergence with an error control. However, this is not standard
even when $d=1$. Unlike \eqref{Newsde}, with \eqref{sde 1.4}, one
needs only to use some usual methods of simulating diffusions as the
ones presented in \cite{glass} for $X$, then simulate $L$ and $R$ as
deterministic functionals of $X$. When both $X$ and the Brownian
motion are $(d\times 1)$-dimensional processes, the contrast between
the two methods become clearer. Add to this the complexity of
studying how the discretization error of \eqref{Newsde} effects the
Longstaff-Schwartz multiple optimal stopping algorithm proposed in
Section \ref{sec42}.

To conclude this section, we would like to point out that the method based on posterior
probabilities is theoretically equivalent to the change of probability method.
Moreover, to solve the problem when $d=1$, one can use some discretization
and weak convergence for both methods (We refer to \cite{KusDup1992} for more details).
Nevertheless, when implementing an algorithm based on Monte carlo as
Longstaff-Schwartz, the use of the change of probability is more appropriate
and could be the method by default when $d > 1$.

\section*{Acknowledgements}
We would like to thank Professor Rama Cont, Professor Paul Feehan,
Professor Steven Shreve, Professor Ivan Yutov, and especially
Professor Huy\^{e}n Pham and Dr. Camelia Pop, for helpful
conversations with the third author. We are grateful to Professor
Hitoshi Ishii for sending to us some hard-to-find classical
literature on viscosity solutions to HJB equations. The first
author's research is supported by MATHEON. The second author's
research is supported by the National Science Foundation under grant
NSF-DMS-09-05754.

\end{document}